\documentclass[a4paper,12pt]{article}
\ifx\TeXtures\undefined
\usepackage[active]{srcltx} 
\usepackage[all]{xy}

\usepackage{hyperref}

\hypersetup{ 
colorlinks=true, 
breaklinks=true, 
urlcolor= black, 
linkcolor= black, 
bookmarksopen=true,
pdftitle={}, pdfauthor={}, pdfsubject={}
}

\usepackage[psamsfonts]{amsfonts}%
\usepackage{amsmath}%
\else
\usepackage{amsmath}%
\usepackage{amsfonts}%
\usepackage{times}
\fi
\usepackage{amsthm} 
\usepackage{epsfig}         
\usepackage[utf8]{inputenc}
\usepackage[T1]{fontenc}
\usepackage[english]{babel}
\usepackage{ae,icomma} 
\usepackage{amssymb,amscd,verbatim}
\usepackage{appendix}
\newcommand{\DATUM}{25-03-2023}              
\newcommand{\change}
{{\marginpar{\#}}}        
\newcommand{\comma}{\: ,}              
\newcommand{\period}{\: .}             
\newcommand{\cA}{{\cal A}}
\newcommand{\cB}{{\cal B}}
\newcommand{\cC}{{\cal C}}

\newcommand{\cI}{{\cal I}}

\newcommand{\cL}{{\cal L}}         
\newcommand{\cP}{{\cal P}}         

\newcommand{\cR}{{\cal R}}
\newcommand{\cS}{{\cal S}}

\newcommand{\cU}{{\cal U}}
\newcommand{\cV}{{\cal V}}
\newcommand{\cW}{{\cal W}}

\newcommand{\field}[1]{\mathbb{#1}}
\newcommand{\R}{\field{R}}            
\newcommand{\Z}{\field{Z}}            
\newcommand{\N}{\field{N}}            
\newcommand{\C}{\field{C}}            

\newcommand{\SSS}{\field{S}}
\newcommand{\uC}{{\underline C}}        
\newcommand{\uP}{{\underline P}}        

\newcommand{\uH}{{\underline H}}

\newcommand{\ux}{{\underline x}}
\newcommand{\uV}{{\underline V}}
\newcommand{\uOmega}{{\underline\Omega}}
\newcommand{\hP}{\widehat{P}}
\newcommand{\hcc}{\hat{c}}
\newcommand{\hu}{\hat{u}}

\newcommand{\tH}{\widetilde{H}}

\newcommand{\tP}{\widetilde{P}}         
\newcommand{\tV}{\widetilde{V}}         

\newcommand{\tc}{\tilde{c}}

\newcommand{\tpsi}{\tilde{\psi}}

\newcommand{\rL}{{\rm L}}                 

\newcommand{\rW}{{\rm W}} 
\newcommand{\rx}{{\rm x}} 
\newcommand{\ry}{{\rm y}}

\newcommand{\cirS}{\mathop{\bigcirc\kern -.73em {\scriptstyle{\rm S}}}}

\newcommand{\cqfd}{\phantom{blablabla}\hfill\qed\newline} 
\newtheorem{theorem}{Theorem}[section]         
\newtheorem{lemma}[theorem]{Lemma}             
\newtheorem{definition}[theorem]{Definition}   
\newtheorem{remark}[theorem]{Remark}           
\newtheorem{proposition}[theorem]{Proposition} 

\theoremstyle{plain}

\newcommand{\donne}{\mapsto}
\newcommand{\dans}{\longrightarrow}

\newcommand{\impl}{\Longrightarrow}

\newcommand{\Pf}{\vspace*{-2mm}{\bf Proof:}\, }
\newcommand{\Pfof}[1]{{\bf Proof of #1:}\, }

\renewcommand{\theequation}{\thesection.\arabic{equation}}

\newcommand{\ncg}{{[\hskip-.7mm [}}
\newcommand{\ncd}{{]\hskip-.7mm ]}}

\begin{document}

\setcounter{section}{0} 

\title{On the analyticity \\
of electronic reduced densities\\
for molecules.}
\author{
{\bf Thierry Jecko}\\
AGM, UMR 8088 du CNRS, site de Saint Martin,\\
2 avenue Adolphe Chauvin,\\
F-95000 Cergy-Pontoise, France. \\
e-mail: jecko@math.cnrs.fr\\
web: http://jecko.perso.math.cnrs.fr/index.html
\\
{\bf Camille Noûs}\\
Laboratoire Cogitamus\\
e-mail: camille.nous@cogitamus.fr\\
web: https://www.cogitamus.fr/
}
\date{\DATUM}
\maketitle

\begin{abstract}
We consider an electronic bound state of the usual, non-relativistic, molecular Hamiltonian with Coulomb interactions and fixed nuclei. Away from appropriate collisions, we prove the real analyticity of all the reduced densities and density matrices, that are associated to this bound state. We provide a similar result for the associated reduced current density. 
\vspace{2mm}

\end{abstract}

\newpage





\section{Introduction.}
\label{intro}
\setcounter{equation}{0}

The notions of density and density matrices are quite old tools in the treatment of physical many body problems. For instance, the Thomas-Fermi theory goes back to the year 1927 (cf. \cite{li}). 
Nowadays these tools do play a central rôle in one important approach of the molecular problem, namely the Density Functional Theory (DFT) (cf. \cite{e,lise}). Regularity properties of the related reduced density and density matrices are useful to develop the DFT. For instance, the real analyticity of a density was recently used for Hohenberg-Kohn Theorems in \cite{g}. The purpose of the present paper is to prove, in a large region of the configuration space, the real analyticity of all reduced densities and density matrices of a pure state for a molecular Hamiltonian with fixed nuclei. 

We consider a molecule with $N$ moving electrons, with $N>1$, and $L$ fixed nuclei, with $L\geq 1$ (Born-Oppenheimer idealization). While the $L$ distinct vectors $R_1, \cdots , R_L\in\R^3$ denote the positions of the nuclei, the positions of the electrons are given by $x_1, \cdots , x_N\in\R^3$. The charges of the nuclei are respectively given by the positive $Z_1, \cdots , Z_L$ and the electronic charge is $-1$. In this picture, the Hamiltonian of the electronic system is 
\begin{eqnarray}\label{eq:hamiltonien}
H&:=&\sum _{j=1}^N\Bigl(-\Delta _{x_j}\, -\, \sum _{k=1}^LZ_k|x_j-R_k|^{-1}\Bigr)\, +\, \sum _{1\leq j<j'\leq N}|x_j-x_{j'}|^{-1}\, +\, E_0\comma\\
\mbox{where}\ E_0&=&\sum _{1\leq k<k'\leq L}Z_kZ_{k'}|R_k-R_{k'}|^{-1}\nonumber
\end{eqnarray}
and $-\Delta _{x_j}$ stands for the Laplacian in the variable $x_j$. Here we denote by $|\cdot|$ the euclidian norm on $\R^3$. Setting 
$\Delta :=\sum _{j=1}^N\Delta _{x_j}$, we define the potential $V$ of the system as the multiplication operator satisfying $H=-\Delta +V$. Thanks to Hardy's inequality 
\begin{equation}\label{eq:hardy}
\exists c>0\, ;\ \forall f\in \rW^{1,2}(\R^3)\comma \ \int _{\R^3}|t|^{-2}\, |f(t)|^2\, dt\ \leq \ c\int _{\R^3}|\nabla f(t)|^2\, dt\comma 
\end{equation}
one can show that $V$ is $\Delta$-bounded with relative bound $0$. Therefore the Hamiltonian $H$ is self-adjoint on the domain of the Laplacian $\Delta $, namely $\rW ^{2,2}(\R^{3N})$ (see Kato's theorem in \cite{rs2}, p. 166-167). \\
From now on, we fix an electronic bound state $\psi\in \rW ^{2,2}(\R^{3N})\setminus\{0\}$ such that, for some real $E$, $H\psi =E\psi$. 
We point out (cf.\ \cite{s,z}) that such a bound state exists at least for appropriate $E\leq E_0$ (cf. \cite{fh}) and for 
\[N\ <\  1\, +\, 2\sum _{k=1}^LZ_k\period\]
Associated to that bound state $\psi$, we consider the following objects. 
Let $k$ be an integer such that $0<k<N$. Let $\rho _k : (\R^3)^k\to\R$ be the almost everywhere defined, $\rL^1(\R^{3k})$-function given by, for $\ux=(x_1; \cdots ; x_k)\in\R^{3k}$, 
\begin{equation}\label{eq:rho-densité}
 \rho _k(\ux)\ =\ \int_{\R^{3(N-k)}}\bigl|\psi (\ux; y)\bigr|^2\, dy\period
\end{equation}
Define $\gamma _k : (\R^3)^{2k}\to\C$ as the almost everywhere defined function given by, for $\ux=(x_1; \cdots ; x_k)\in\R^{3k}$ and $\ux'=(x_1'; \cdots ; x_k')\in\R^{3k}$, 
\begin{equation}\label{eq:gamma-densité}
 \gamma _k(\ux; \ux')\ =\ \int_{\R^{3(N-k)}}\overline{\psi (\ux; y)}\psi (\ux'; y)\, dy\period
\end{equation}
By \cite{k}, $\psi$ is actually a continuous function. Therefore, the objects $\rho _k$ and $\gamma _k$ are everywhere defined and, $\rho _k(\ux)=\gamma _k(\ux; \ux)$ holds true everywhere. 
It turns out that $\gamma _k$ may be seen as the kernel (``Green function'') of a trace class operator (see \cite{le,lise}). We call the function $\rho _k$ the {\em $k$-particle reduced density} and the kernel $\gamma _k$ the {\em $k$-particle reduced density matrix}. \\
When $k=1$, $\rho:=\rho_1$ is often simply called the (electronic) density. In this case, 
we also introduce the {\em reduced current density}, defined as the almost everywhere defined, $\rL^1(\R^{3})$-function $C : \R^3\to\R^3$ given by, for $x\in\R^{3}$, the following imaginary part 
\begin{equation}\label{eq:courant-densité}
 C(x)\ =\ \Im\, \int_{\R^{3(N-1)}}\overline{\nabla_x\psi (x; y)}\psi (x; y)\, dy\period
\end{equation}
From a physical point of view, the previous objects differ from the true physical ones by some prefactor (see \cite{e,le,lise,lsc}). \\
It is useful to introduce the following subsets of $\R^{3k}$. Denoting for a positive integer $p$ by $\ncg 1; p\ncd$ the set of the integers $j$ satisfying $1\leq j\leq k$, the closed set 
\[\cC _k\ :=\ \bigl\{\ux=(x_1; \cdots ; x_k)\in\R^{3k}\, ;\, \exists (j; j')\in\ncg 1; k\ncd ^2\, ;\, j\neq j'\ \mbox{and}\ x_j=x_{j'}\bigr\}\]
gathers all possible collisions between the first $k$ electrons while the closed set 
\[\cR _k\ :=\ \bigl\{\ux=(x_1; \cdots ; x_k)\in\R^{3k}\, ;\, \exists j\in\ncg 1; k\ncd \, ,\, \exists \ell\in\ncg 1; L\ncd \, ;\, x_j=R_\ell\bigr\}\]
groups together all possible collisions of these $k$ electrons with the nuclei. We set $\cU^{(1)}_k:=\R^{3k}\setminus (\cC_k\cup\cR_k)$, which is an open subset of $\R^{3k}$. \\
Now, we consider two sets of positions for the first $k$ electrons and introduce the set of all possible collisions between positions of differents sets, namely  
\[\cC _k^{(2)}\ :=\ \left\{
\begin{array}{l}
 (\ux ; \ux ')\in(\R^{3k})^2\, ;\ \ux=(x_1; \cdots ; x_k)\, ,\ \ux '=(x_1'; \cdots ; x_k')\, ,\, \\
 \\
 \exists (j; j')\in\ncg 1; k\ncd ^2\, ;\  x_j=x_{j'}'
\end{array}
\right\}\period
\]
We introduce the open subset of $(\R^{3k})^2$ defined by 
\[\cU^{(2)}_k\ =\ \bigl(\cU^{(1)}_k\times\cU^{(1)}_k\bigr)\setminus\cC _k^{(2)}\period\]
In the $k=1$ case, we note that $\cC_1=\emptyset$, $\cU^{(1)}_1=\R^3\setminus\{R_1,\cdots , R_L\}$ and 
\[\cU^{(2)}_1\ =\ \bigl(\R^3\setminus\{R_1,\cdots , R_L\}\bigr)^2\setminus D\comma\]
where $D=\{(x; x')\in(\R^{3})^2; x=x'\}$ is the diagonal of $(\R^{3})^2$. 
We have the following
\begin{theorem}\label{th:rho-anal}{\rm \cite{fhhs1,fhhs2,j}}.\\
The one-particle reduced density $\rho :=\rho _1$ is real analytic on $\cU_1^{(1)}=\R^3\setminus\{R_1,\cdots , R_L\}$.
\end{theorem}
The original proof was given in \cite{fhhs1,fhhs2} and relies on a clever 
decomposition of the integration domain of the integral defining $\rho$ into pieces, on which one can perform an appropriate change of variables. A different proof of Theorem~\ref{th:rho-anal} is provided in \cite{j}. We refer to \cite{fhhs3,fs} for a refined information on $\rho$. 
In \cite{fhhs2}, a similar analyticity result on $\gamma _1$ was announced, without proof. More precisely, it was claimed there that $\gamma _1$ should be real analytic on $\cU^{(1)}_1\times\cU^{(1)}_1$. Adapting the method of \cite{fhhs2}, a shlighly weaker result on $\gamma_1$ was proved in \cite{hs}, namely 
\begin{theorem}\label{th:gamma-1-anal}{\rm \cite{hs}}.\\
The one-particle reduced density matrix $\gamma _1$ is real analytic on $\cU_1^{(2)}$.
\end{theorem}
Note that $\cU_1^{(2)}$ is indeed a strict subset of $\, \cU^{(1)}_1\times\cU^{(1)}_1$. \\
In \cite{fhhs2}, it was also announced without proof that $\rho _2$ should be real analytic on 
\[\cU^{(1)}_2\ =\ \R^{6}\setminus (\cC_2\cup\cR_2)\ =\ \bigl(\R^3\setminus\{R_1,\cdots , R_L\}\bigr)^2\setminus D\ =\ \cU^{(2)}_1\period\]
In the present paper, we work with the tools and the method, used in \cite{j} to get Theorem~\ref{th:rho-anal}, to prove the following three results. 
\begin{theorem}\label{th:courant-anal}
The reduced current density $C$ is real analytic on $\cU_1^{(1)}=\R^3\setminus\{R_1,\cdots , R_L\}$.
\end{theorem}
\begin{theorem}\label{th:rho-k-anal}
For all integer $k$ with $0<k<N$, the $k$-particle reduced density $\rho _k$ is real analytic on $\cU^{(1)}_k$. 
\end{theorem}
\begin{theorem}\label{th:gamma-k-anal}
 For all integer $k$ with $0<k<N$, the $k$-particle reduced density matrix $\gamma _k$ is real analytic on $\cU^{(2)}_k=(\cU^{(1)}_k\times\cU^{(1)}_k)\setminus\cC _k^{(2)}$. 
\end{theorem}
Theorem~\ref{th:rho-k-anal} for $k=1$ coincides with Theorem~\ref{th:rho-anal}. The proof of Theorem~\ref{th:rho-k-anal} below is, when restricted to the $k=1$ case, just a rewriting of the proof of Theorem~\ref{th:rho-anal} in \cite{j} and provides a direct justification of Theorem~\ref{th:courant-anal}. 
We note that Theorem~\ref{th:gamma-k-anal} for $k=1$ is exactly Theorem~\ref{th:gamma-1-anal} and that Theorem~\ref{th:rho-k-anal} with $k=2$ fits precisely to the announced result on $\rho _2$ in \cite{fhhs2}. \\
The set of all possible collisions between particles is $\cC_N\cup\cR_N$ and the potential $V$ is real analytic precisely on $\R^{3N}\setminus (\cC_N\cup\cR_N)$. Classical elliptic regularity applied to the equation $H\psi=E\psi$ shows that $\psi$ is also real analytic on $\R^{3N}\setminus (\cC_N\cup\cR_N)$. We refer to \cite{acn,fhhs3,fhhs4,fhhs5,fs,k} for more information on the behaviour of $\psi$ near collisions. It turns out that, at least at some places in $\cC_N\cup\cR_N$, $\psi$ is not real analytic (cf. \cite{fhhs3,k}). Thus, in the definition \eqref{eq:rho-densité} of $\rho _k$, for an integer $k$ with $0<k<N$, the domain of integration does contain such singularities. Therefore it is not clear a priori that $\rho _k$ is real analytic somewhere. This is however true, by Theorem~\ref{th:rho-k-anal}, on the complement of the set of all possible collisions of the ``external'' electrons, namely on $\cU^{(1)}_k$. \\
The proofs of Theorems~\ref{th:rho-anal} and \ref{th:gamma-1-anal} do not directly use the analyticity of $\psi$. This is also the case of the proofs of Theorems~\ref{th:courant-anal}, \ref{th:rho-k-anal}, and~\ref{th:gamma-k-anal}, below. \\
It turns out that our proofs of Theorems~\ref{th:rho-k-anal} and~\ref{th:gamma-k-anal} have a common structure that we want to describe now. Using an appropriate, $\ux$-dependent, unitary operator $U_\ux$, that acts on $\rL^2(\R^{3(N-k)})$, we locally transform the equation $H\psi=E\psi$ into an elliptic differential equation $(U_\ux HU_\ux^{-1})U_\ux\psi=EU_\ux\psi$ (the ``twisted'' equation), the coefficients of which have nice analytic properties. Applying elliptic regularity to the latter equation, we obtain the real analyticity of $\ux\donne U_\ux\psi\in\rL^2(\R^{3(N-k)})$ that yields the real analyticity of $\rho_k : \ux\donne\|U_\ux\psi\|^2$ ($\|\cdot\|$ being the norm on $\rL^2(\R^{3(N-k)}))$. For $\gamma_k$, we follow the same lines but with a $(\ux ; \ux ')$-dependent, unitary operator $U_{(\ux ; \ux ')}$, and two ``twisted'' differential equations in the variables $(\ux ; \ux '; y)$. These proofs are performed in Section~\ref{s:proof}. \\
Let us now compare the two available methods to prove such results, namely the one mentioned just above and the one used in \cite{fhhs2,hs} to get Theorems~\ref{th:rho-anal} and \ref{th:gamma-1-anal}. We point out that elliptic regularity was also an important tool in \cite{fhhs2,hs}. Another common tool is the use of a $\ux$-dependent change of variables on the ``internal'' variables $y\in \R^{3(N-k)}$, since the above unitary operators are both the implementation in $\rL^2(\R^{3(N-k)})$ of such a change of variables. However, we use here global changes of variables in contrast to the local ones in \cite{fhhs2,hs}. While, there, the integrals on $\R^{3(N-1)}$ defining $\rho _1$ and $\gamma _1$ are split in several pieces, on which an appropriate change of variables is performed to make the real analyticity apparent, we view here $\rho _k$ (resp. $\gamma _k$) as the $\rL^2(\R^{3(N-k)})$-norm (resp. the $\rL^2(\R^{3(N-k)})$-scalar product) of one (resp. two) $\rL^2(\R^{3(N-k)})$-valued, real analytic function(s). \\
In view of the behaviour of $\psi$ near some points in $\cC_N\cup\cR_N$, that was established in \cite{fhhs5}, it is quite natural to expect that, for an integer $k$ with $0<k<N$, neither $\rho _k$ nor $\gamma _k$ is real analytic everywhere. Their exact domain of real analyticity is still an open question. The domains obtained in Theorems~\ref{th:rho-k-anal} and~\ref{th:gamma-k-anal} are the largest ones that can be reached by the method used in this paper, as explained in Remarks~\ref{r:limitation-1} and~\ref{r:limitation-2}. 

As already pointed out in \cite{j}, the twisted calculus, that we use here, allows us to treat more general Hamiltonians than $H$. We introduce, in Section~\ref{s:extensions}, a class of elliptic, second order differential operators with Coulomb singularities. We observe that the conjugation by the twist actually preserves the ellipticity. Assuming that there exists a bound state for such an operator, we show in Theorem~\ref{th:géné-anal} that our previous, main results hold true for it. 
Moreover, we expect that this should extend to a class of elliptic, {\em pseudo}-differential operators with Coulomb singularities, since the used tools are available for such operators. We refer to \cite{ms} for more information on the twisted (pseudo-)differential calculus. 

Maybe it would be useful for the DFT to extend the present results on pure states to general states (see Sections 3.1.4 and 3.1.5 in \cite{lise} for details). We expect that Theorem~\ref{th:regu-twisted-op} below could be a good starting point for this purpose. 

The paper is organized as follow. In Section~\ref{s:elliptic}, we provide a general notation, basic facts on real analyticity, and a result on elliptic regularity from \cite{j}. In Section~\ref{s:twisted-diff}, we review the twisted calculus based on Hunziker's twist, show how it can be used to ``desingularize'' a differential equation with Coulomb singularities, and apply elliptic regularity to the corresponding twisted equation. Section~\ref{s:proof} is devoted to the proof of our main results, namely Theorems~\ref{th:courant-anal}, \ref{th:rho-k-anal}, and~\ref{th:gamma-k-anal}. The extension of these results to a larger class of Hamiltonians is performed in Section~\ref{s:extensions}. Finally, we gathered in Appendix~\ref{app:computations} some elementary computations that are used in the main text and provide in Appendix~\ref{app:pseudo} basic material on global pseudo-differential operators. 

{\bf Acknowledgments:} The authors thanks S. Fournais and T. \O stergaard S\o rensen for fruitful discussions related to the subject of this paper. Dedicated to K. 

\section{Basic tools and elliptic regularity.}
\label{s:elliptic}
\setcounter{equation}{0}

To prepare the proof of the main results, we recall basic tools and state a general result on the analytic regularity of solutions of some elliptic equations.

\subsection{Notation and basic tools.}
\label{ss:notation}

We start with a general notation. We denote by $\C$ the field of complex numbers, and, for $z\in\C$, $\Re z$ (resp. $\Im z$) stands for the real part (resp. imaginary part) of $z$. \\
Let $p$ be a positive integer. Recall that, for $u\in\R^p$, we write $|u|$ for the euclidian norm of $u$. Given such a vector $u\in\R^p$ and a nonnegative real number $r$, we denote by $B(u; r[$ (resp. $B(u; r]$) the open (resp. closed) ball of radius $r$ and centre $u$, for the euclidian norm $|\cdot|$ in $\R^p$. \\
In the one dimensional case, we use the following convention for (possibly empty) intervals: for $(a; b)\in\R^2$, let $[a; b]=\{t\in\R; a\leq t\leq b\}$, $[a; b[=\{t\in\R; a\leq t<b\}$, $]a; b]=\{t\in\R; a<t\leq b\}$, and $]a; b[=\{t\in\R; a<t<b\}$.\\
We denote by $\N$ the set of nonnegative integers and set $\N^\ast=\N\setminus\{0\}$. 
If $p\leq q$ are non negative integers, we set $\ncg p; q\ncd :=[p; q]\cap\N$, $\ncg p; q\ncg =[p; q[\cap\N$, $\ncd p; q\ncg =]p; q[\cap\N$, and $\ncd p; q\ncd :=]p; q]\cap\N$.\\
Given an open subset $O$ of $\R^p$ and $n\in\N$, we denote by $W^{n,2}(O)$ the standard Sobolev space of those $\rL^2$-functions on $O$ such that, for $n'\in \ncg 0; n\ncd$, their distributional derivatives of order $n'$ belong to $\rL^2(O)$. We denote by $\langle\cdot , \cdot\rangle _n$ (resp. $\|\cdot\|_n$) the right linear scalar product (resp. the norm) on $W^{n,2}(O)$. In particular, $W^{0,2}(O)=\rL^2(O)$. Without reference to $p$ and $O$, we denote by $\|\cdot\|$ (resp. $\langle\cdot , \cdot\rangle$) the $\rL^2$-norm (resp. the right linear scalar product) on $\rL^2(O)$. \\
For two Banach spaces $(\cA , \|\cdot\|_\cA)$ and $(\cB , \|\cdot\|_\cB)$, the space $\cL (\cA ; \cB)$ of continuous linear maps from $\cA$ to $\cB$ is also a Banach space for the operator norm $\|\cdot\|_{\cL (\cA ;\, \cB)}$ defined by 
\[\forall M\in\cL (\cA ; \cB)\comma\hspace{.4cm}\|M\|_{\cL (\cA ;\, \cB)}\ :=\ \sup_{a\in\cA\setminus\{0\}}\, \frac{\|M(a)\|_\cB}{\|a\|_\cA}\period\]
We simply denote $\cL (\cA ;\, \cA)$ by $\cL (\cA)$. \\
Let $p$ be a positive integer and $O$ an open subset of $\R^p$. Let $(\cA , \|\cdot\|_\cA)$ be a Banach space and $\theta : O\ni\rx=(\rx_1; \cdots ; \rx _p)\donne\theta (\rx)\in\cA$. For $j\in\ncg1; p\ncd$, we denote by $\partial _j\theta$ or $\partial _{\rx_j}\theta$ the $j$'th first partial derivative of $\theta$. For $\alpha\in\N^p$ and $\rx\in\R^d$, we set $D_\rx^\alpha :=(-i\partial _\rx)^\alpha :=(-i\partial _{\rx_1})^{\alpha _1}\cdots(-i\partial _{\rx_p})^{\alpha _p}$, $D_\rx=-i\nabla_\rx$, $\rx^\alpha :=\rx_1^{\alpha_1}\cdots \rx_p^{\alpha_p}$, $|\alpha |:=\alpha_1+\cdots +\alpha _p$, $\alpha !:=(\alpha_1!)\cdots (\alpha_p !)$, $|\rx|^2=\rx_1^2+\cdots +\rx_p^2$, and $\langle \rx\rangle :=(1+|\rx|^2)^{1/2}$. \\
We choose the same notation for the length $|\alpha|$ of a multiindex $\alpha\in\N^p$ and for the euclidian norm $|\rx|$ of a vector $\rx\in\R^p$ but the context should avoid any confusion. \\
We denote by $\cC^\infty(O;\cA)$ (resp. $\cC_b^\infty(O;\cA)$, resp. $\cC_c^\infty(O;\cA)$, resp. $\cC^\omega(O;\cA)$) the vector space of functions from $O$ to $\cA$ which are smooth (resp. smooth with bounded derivatives, resp. smooth with compact support, resp. real analytic). In the case $\cA=\C$, we simply write $\cC^\infty(O)$ (resp. $\cC_b^\infty(O)$, resp. $\cC_c^\infty(O)$, resp. $\cC^\omega(O)$). 

Now we recall basic facts on analytic functions. Take again a Banach space $(\cA , \|\cdot\|_\cA)$, a positive integer $p$, and an open subset $O$ of $\R^p$. For a smooth function $\theta : O\dans\cA$, the following properties are equivalent: $\theta \in\cC^\omega(O;\cA)$; 
\begin{equation}\label{eq:anal-1}
\exists A>0\comma\, \forall\, \mbox{compact}\, K\subset O\comma\, \forall \alpha\in\N^p\comma \hspace{.4cm} \sup _{\rx\in K}\, 
\bigl\|(D_\rx^\alpha \theta )(\rx)\bigr\|_{\cA}\ \leq \ A^{|\alpha |+1}\cdot (\alpha !)\, ;
\end{equation}
\begin{equation}\label{eq:anal-2}
 \exists A>0\comma\, \forall\, \mbox{compact}\, K\subset O\comma\, \forall \alpha\in\N^p\comma \hspace{.4cm} \sup _{\rx\in K}\, 
\bigl\|(D_\rx^\alpha \theta )(\rx)\bigr\|_{\cA}\ \leq \ A^{|\alpha |+1}\cdot (|\alpha | !)\, ;
\end{equation}
\begin{equation}\label{eq:anal-3}
 \exists A>0\comma\, \forall\, \mbox{compact}\, K\subset O\comma\, \forall \alpha\in\N^p\comma \hspace{.4cm} \sup _{\rx\in K}\, 
\bigl\|(D_\rx^\alpha \theta )(\rx)\bigr\|_{\cA}\ \leq \ A^{|\alpha |+1}\cdot (1+|\alpha |)^{|\alpha |}\period
\end{equation}

We refer to \cite{h3,hs,j} for details. \\
Let $(\cB , \|\cdot\|_\cB)$ be another Banach space and $f : \cA^2\dans\cB$ be a continuous bilinear map. Then, for any $(\theta_1 ; \theta _2)\in (\cC^\omega(O;\cA))^2$, the map $f(\theta_1 ; \theta _2): O\ni\rx\donne f(\theta_1 (\rx); \theta _2(\rx))\in\cB$ is real analytic. \\
Consider the case where $\cA=W^{n,2}(O)$, for some integer $n$, and recall that $\langle\cdot , \cdot\rangle _n$ (resp. $\|\cdot\|_n$) is the right linear scalar product (resp. the norm) on $\cA$. Take $(\theta_1 ; \theta _2)\in (\cC^\omega(O;\cA))^2$. Then the map $\langle\theta _1, \theta _2\rangle_n : O\dans\C$, defined by $\langle\theta _1, \theta _2\rangle_n (\rx)=\langle\theta _1(\rx), \theta _2(\rx)\rangle_n$, is real analytic. So is also the real valued map $\|\theta _1\|_n^2 : O\ni\rx\donne \|\theta _1(\rx)\|_n^2$. In particular, taking $n=0$, the maps $\langle\theta _1, \theta _2\rangle : O\ni\rx\donne \langle\theta _1(\rx), \theta _2(\rx)\rangle\in\C$ and $\|\theta _1\|^2 : O\ni\rx\donne \|\theta _1(\rx)\|^2\in\R^+$ are also real analytic. 

Finally we need some specific notation to describe the structure of the considered Hamiltonian and of our results.\\
For a $\cA$-valued function $\theta : \R^3\ni z=(z^{(1)}; z^{(2)}; z^{(3)})\donne\theta (z)\in\cA$, we denote by $d\theta$ or $d_z\theta$ the total derivative of $\theta$. According to the previous notation, we denote, for $j\in\{1; 2; 3\}$, by $\partial _j\theta$ or $\partial _{z^{(j)}}\theta$ the $j$'th first partial derivative. For a multiindex $\alpha =(\alpha _1; \alpha _2; \alpha _3)\in\N^3$, we set $D_z^\alpha :=(-i\partial _z)^\alpha :=(-i\partial _1)^{\alpha _1}(-i\partial _2)^{\alpha _2}(-i\partial _3)^{\alpha _3}$, $D_z=-i\nabla _z$, $z^\alpha=(z^{(1)})^{\alpha _1}(z^{(2)})^{\alpha _2}(z^{(3)})^{\alpha _3}$, and $|\alpha|=\alpha _1+\alpha _2+\alpha _3$.  \\
Let $m$ be a positive integer. We write a point $\rx\in\R^{3m}$ as $\rx=(\rx _1; \cdots ; \rx_ m)$. For any multiindex $\alpha =(\alpha _1; \cdots ; \alpha _m)\in\N^{3m}$, we set $D_\rx^\alpha :=(-i\partial _\rx)^\alpha :=(-i\partial _{\rx _1})^{\alpha _1}\cdots (-i\partial _{\rx _m})^{\alpha _m}$, $D_\rx=-i\nabla _\rx$, $\rx^\alpha=\rx_1^{\alpha _1}\cdots\rx_m^{\alpha _m}$, and $|\alpha |=|\alpha _1|+\cdots+|\alpha_m|$.\\
Let $m$ and $p$ be positive integers. We write $\rx\in\R^{3m}$ (resp. $\ry\in\R^{3p}$) as $\rx=(\rx _1; \cdots ; \rx_ m)$ (resp. $\ry=(\ry _1; \cdots ; \ry_ p)$). 
For a $\cA$-valued function $\theta :\R^{3m}\times\R^{3p}\ni (\rx; \ry)\donne \theta (\rx; \ry)\in\cA$, let $d_\rx\theta$ (resp. $d_\ry \theta$) be the total derivative of $\theta$ w.r.t.\ $\rx$ (resp. $\ry$). \\
For $n\in\N$, $W^{n,2}(\R^{3p})$ denotes the standard Sobolev space of $\rL^2$-functions on $\R^{3p}$ admitting, up to order $n$, $\rL^2$ distributional derivatives. In particular, $W^{0,2}(\R^{3p})=\rL^2(\R^{3p})$. We set $\cW_n=W^{n,2}(\R^{3p})$ and let $\cB_n=\cL (\cW_n; \cW_0)$ be the Banach space of linear continuous maps from $\cW_n$ to $\cW_0$.

\subsection{Elliptic regularity.}
\label{ss:elliptic-regu}

In this subsection, we extend, essentially with the same proof, a regularity result that was proved but not stated in \cite{j}. 

Let $(m; p)\in (\N^\ast)^2$. Let $\Omega$ be an open subset of $\R^{3m}$. We consider second order differential operators $P$ on $\Omega\times\R^{3p}$ given by 
\begin{equation}\label{eq:op-diff2}
 P\ =\ \sum_{\stackrel{(\alpha ; \, \beta)\in\N^{3m}\times\N^{3p}}{\stackrel{}{|\alpha |+|\beta |\leq 2}}}\, c_{\alpha\beta}(\rx; \ry)\, D_\rx^\alpha D_\ry^\beta\comma
\end{equation}
where the coefficients $c_{\alpha\beta}$ belong to $\cC_b^\infty(\Omega\times\R^{3p}; \C)$. 
In particular, for $\rx\in\Omega$, the multiplication operator $C_{\alpha\beta}(\rx ; \cdot)$ by $c_{\alpha\beta}(\rx; \cdot)$ belongs to $\cL (\cW _n)$, for all $n\in\N$. \\
The principal symbol of those operators $P$ is the smooth, complex-valued function $\sigma _P$, that is defined on $\Omega\times\R^{3p}\times\R^{3m}\times\R^{3p}$ by 
\begin{equation}\label{eq:symb-principal}
\sigma _P(\rx ; \ry ; \xi ; \eta)\ =\ \sum_{\stackrel{(\alpha ;\, \beta)\in\N^m\times\N^p}{|\alpha |+|\beta |=2}}\, c_{\alpha\beta}(\rx; \ry)\, \xi^\alpha \eta^\beta\period
\end{equation}
We assume that $P$ is (globally) elliptic in the following sense: there exists $C_P>0$ such that, for all $(\rx ; \ry ; \xi ; \eta)\in\Omega\times\R^{3p}\times\R^{3m}\times\R^{3p}$, 
\begin{equation}\label{eq:glob-ellipt}
 |\xi|^2\, +\, |\eta|^2\, \geq \, 1\hspace{.4cm}\impl\hspace{.4cm}\bigl|\Re\sigma _P(\rx ; \ry ; \xi ; \eta)\bigl|\ \geq\ C_P\bigl(|\xi|^2\, +\, |\eta|^2\bigr)\period
\end{equation}
Furthermore, we require that, for $(\alpha ; \beta)\in\N^{3m}\times\N^{3p}$ with 
$|\alpha |+|\beta |\leq 2$, the map $\Omega\ni\rx\donne C_{\alpha\beta}(\rx ; \cdot)\in\cB_0$ is well-defined and real analytic.\\
Under these assumptions on $P$, we note that, for all $\alpha\in\N^{3m}$ with $|\alpha |\leq 2$, the function $a_\alpha : \Omega\dans\cB_{2-|\alpha |}$, which maps each $\rx\in\Omega$ to the differential operator 
\begin{equation}\label{eq:coeff-op}
a_\alpha(\rx)\ =\ \sum_{\stackrel{\beta\in\N^{3p}}{\stackrel{}{|\beta |\leq 2-|\alpha |}}}\, c_{\alpha\beta}(\rx; \cdot)\, D_\ry^\beta
\end{equation}
on $\R^{3p}$, is well-defined and real analytic. 
\begin{definition}\label{def:ope-diff}
 For $(m; p)\in (\N^\ast)^2$ and $\Omega$ an open subset of $\R^{3m}$. We denote by ${\rm Diff}_2(\Omega)$ the set of all the differential operators $P$ on $\Omega\times\R^{3p}$ satisfying the previous requirements.
\end{definition}
\begin{theorem}\label{th:ellipti-regu}{\rm\cite{j}}.\\
Let $P$ be differential operator in the class ${\rm Diff}_2(\Omega)$. Let $W : \Omega\dans\cB_1$ be a real analytic map. Take $\varphi\in W^{2,2}(\Omega\times\R^{3p})$ such that $(P+W)\varphi =0$. Then, the map $\Omega\ni\rx\donne\varphi (\rx; \cdot)\in\cW_2$ is real analytic. 
\end{theorem}
\begin{remark}\label{r:elliptique-class}
 If, in Theorem~\ref{th:ellipti-regu}, one replaces $3m$ by some positive integer $n$ and one sets $p=0$ and $W=0$ with the convention that $\Omega\times\R^0=\Omega$ and $\cW_2=\C$, then one recovers exactly Hörmander's Theorem 7.5.1 in {\rm\cite{h1}}. We assign Theorem~\ref{th:ellipti-regu} to {\rm\cite{j}} since, although it was not stated there, it was proved in a particular case and the corresponding proof directly extends to the present, general case. In {\rm \cite{fhhs1,fhhs2,hs}}, one can find similar results in a slightly narrower framework. Therefore, Theorem~\ref{th:ellipti-regu} is more or less known in the literature. 
\end{remark}
Nevertheless, we provide below, for completeness, a proof of Theorem~\ref{th:ellipti-regu}. To make it more transparent and facilitate the comparison with known results of this kind, we state intermediate results in some lemmata. \\
For convenience, we choose to use an appropriate pseudo-differential calculus. Since we are concerned with a local result (in $\Omega$), it is natural to take a local one (see, for instance, in \cite{h2} on page 83 or in \cite{t}, section 7.10). However, we prefer to use the basic global calculus presented in the beginning of Chapter 18 in \cite{h2}, since we think that it is more accessible. \\
Let $P\in{\rm Diff}_2(\Omega)$. To prove Theorem~\ref{th:ellipti-regu}, it suffices to show the required real analyticity near each point in $\Omega$. Let $\rx^0\in\Omega$ and let $\Omega _0'$ be a bounded neighbourhood of $\rx^0$ such that $\Omega _0'\subset\Omega$. As a first step, we construct an extension $\hP\in{\rm Diff}_2(\R^{3m})$ of the restriction of $P$ to $\Omega _0'$. 
\\
In view of \eqref{eq:glob-ellipt}, we observe that the real part $\Re\sigma_P$ of $\sigma_P$ does not vanish, except for $(\xi ; \eta)=(0;0)$. Since it is continuous, it must keep the same sign $\sigma^0$ on 
\[\Omega\times\R^{3p}\times (\R^{3(m+p)}\setminus\{(0; 0)\})\period\]
Furthermore, \eqref{eq:glob-ellipt} also implies that, for $(\rx ; \ry; \xi ; \eta)\in\Omega\times\R^{3p}\times\R^{3(m+p)}$, 
\begin{equation}\label{eq:glob-ellipt-2}
\sigma^0\cdot\Re\sigma _P(\rx ; \ry ; \xi ; \eta)\ \geq\ C_P\bigl(|\xi|^2\, +\, |\eta|^2\bigr)\period
\end{equation}
Let $\chi\in\cC_c^\infty(\R^{3m})$ be a cut-off function with support in $\Omega$ such that $0\leq\chi\leq 1$ and, on $\Omega _0'$, $\chi =1$. Recalling that $P$ is given by \eqref{eq:op-diff2}, we set 
\[\hP\ =\ \sum_{\stackrel{(\alpha ; \, \beta)\in\N^{3m}\times\N^{3p}}{\stackrel{}{|\alpha |+|\beta |\leq 2}}}\, \hcc_{\alpha\beta}(\rx; \ry)\, D_\rx^\alpha D_\ry^\beta\comma\]
where $\hcc_{\alpha\beta}(\rx; \ry)=\chi (\rx)c_{\alpha\beta}(\rx; \ry)$, for all $(\rx ; \ry)\in\R^{3(m+p)}$ and $|\alpha|+|\beta|<2$, and where, for all $(\rx ; \ry ; \xi ; \eta)\in\R^{6(m+p)}$, 
\begin{align}
\sigma _{\hP} (\rx ; \ry ; \xi ; \eta)\ &:=\ \sum_{\stackrel{(\alpha ;\, \beta)\in\N^m\times\N^p}{|\alpha |+|\beta |=2}}\, \hcc_{\alpha\beta}(\rx; \ry)\, \xi^\alpha \eta^\beta\nonumber\\ 
\ &:=\ \chi (\rx)\cdot\sigma _P (\rx ; \ry ; \xi ; \eta)\, +\, \sigma^0C_P\cdot \bigl(1\, -\, \chi (\rx)\bigr)\cdot \bigl(|\xi|^2\, +\, |\eta|^2\bigr)\period\label{eq:symb-principal-2}
\end{align}
We observe that \eqref{eq:glob-ellipt-2} with $\sigma _P$ replaced by $\sigma _{\hP}$ holds true on $\R^{6(m+p)}$. In particular, $\hP$ is elliptic. We immediately verify that $\hP\in{\rm Diff}_2(\R^{3m})$. We note that $\hP$ coincides with $P$ as differential operator on $\Omega_0'$. \\
We observe that the symbol of $\hP$ belongs to the global class $S^2(\R^{6N})$ of symbols, for $N=m+p$, and, since it is elliptic, we can construct a parametrix (cf. Appendix~\ref{app:pseudo}). Thus there exist two pseudo-differential operators $Q=q(\rx; \ry; D_\rx ; D_\ry)$ and $R=r(\rx; \ry; D_\rx ; D_\ry)$ 
with symbols $q, r\in S^{-2}(\R^{6N})$ such that $Q\hP=I-R$ (where $I$ is the identity differential operator) and 
\begin{equation}\label{eq:regularisation}
Q\comma \, R\, \in\,  \cL\bigl(\rW^{k, 2}(\R^{3N}); \rW^{k+2, 2}(\R^{3N})\bigr)\period
\end{equation}
\Pfof{Theorem~\ref{th:ellipti-regu}} We take $\varphi\in W^{2,2}(\Omega\times\R^{3p})$ such that $(P+W)\varphi =0$. Let $\rx^0\in\Omega$. We use the above construction. Let $\Omega_0$ a neighbourhood of $x^0$ such that $\Omega_0\subset\Omega_0'$. Let $\chi _0\in C_c^\infty(\R^3)$ with $\chi _0=1$ near $\Omega_0$ and with support in $\Omega _0'$.  Thus $\chi\chi_0=\chi_0$, where $\chi$ is the cut-off function appearing in the definition of $\hP$.
\begin{lemma}\label{l:smooth}
We have $\chi_0\varphi\in\cC_c^\infty(\R^{3m}; \cW_2)$.  
\end{lemma}
\Pf  Since $Q\hP=I-R$, $\hP\chi_0=\chi P\chi_0$, and $(P+W)\varphi=0$, we obtain 
\begin{equation}\label{elliptique-chi_0}
\chi _0\varphi \ =\ R\chi _0\varphi\, +\, Q\chi [P\, ,\, \chi _0]\varphi \, -\, Q\chi _0W\varphi\comma
\end{equation}
where the commutator $[P\, ,\, \chi _0]$ is a first order diffential operator in $\rx$. 
Since $\varphi\in W^{2,2}(\Omega\times\R^{3p})$, $\varphi\in W^{1,2}(\Omega ; \cW_1)$ and 
$\chi_0\varphi\in W^{1,2}(\R^{3m}; \cW _1)$. 
We successively observe that 
\[W\varphi\in W^{1,2}(\Omega; \cW _0)\comma\ \chi _0W\varphi\in W^{1,2}(\R^{3m}; \cW _0)\comma\hspace{.4cm}\mbox{and}\hspace{.4cm} Q\chi W\chi_0\varphi\in W^{2,2}(\R^{3m}; \cW _1)\period\]
Similarly, we have 
\[[P\, ,\, \chi _0]\varphi\in W^{0,2}(\Omega; \cW _1)\comma\ \chi [P\, ,\, \chi _0]\varphi\in W^{0,2}(\R^{3m}; \cW _1)\comma\]
\[Q\chi [P\, ,\, \chi _0]\varphi\in W^{2,2}(\R^{3m}; \cW _1)\hspace{.4cm}\mbox{and}\hspace{.4cm}R\chi_0\varphi\in W^{3,2}(\R^{3m}; \cW _1)\period\]
Thus $\chi_0\varphi\in W^{2,2}(\R^{3m}; \cW _1)$, by \eqref{elliptique-chi_0}. Now, we apply the same argument again but with this new information, yielding $\chi_0\varphi\in W^{2,2}(\R^{3m}; \cW _2)$. In this way, we prove by induction that, for all $k\in\N$, $\chi_0\varphi\in W^{k,2}(\R^{3m}; \cW _2)$, yielding $\chi_0\varphi\in\cC^\infty(\R^{3m}; \cW_2)$. \cqfd
\begin{lemma}\label{l:bound-derivative}
 There exists $C_e>0$ such that, for all $v\in \cC _c^\infty (\R^{3m};\cW _2)$ with support in $\Omega _0'$, for all $r\in\{0;1;2\}$, for all $\alpha\in\N^{3m}$,
\begin{equation}\label{eq:borne-elliptique}
|\alpha |+r\leq 2\ \impl \ \bigl\|D_\rx^\alpha v\bigr\|_{\rL ^2(\Omega _0;\cW _r)}\ \leq \ C_e\bigl\|(P+W)v\bigr\|_{\rL ^2(\Omega _0;\cW _0)}\period
\end{equation}
\end{lemma}
\Pf Let $P_2:=\sigma _P(\rx ; \ry ; D_\rx ; D_\ry)$ and $P_1=P-P_2$. We consider the following operators as operators in $\rL ^2(\Omega _0;\cW _0)$ with domain $W^{2,2}(\Omega _0\times\R^{3p})$. By ellipticity, $\Delta=\Delta _\rx+\Delta _\ry$ is $P_2$-bounded. By assumption, $W$ is $(-\Delta_\rx)^{1/2}$-bounded. 
Since the operators $(-\Delta_\rx)^{1/2}$ and $D_\rx^\alpha D_\ry^\beta$, for $|\alpha|+|\beta|\leq 1$, are $\Delta$-bounded with relative bound less than one, 
so are also $W$ and, thanks to the boundedness assumption on its coefficients, $P_1$. Thus $W$ and $P_1$ are $P_2$-bounded with relative bound less than one. This implies that $W$ is $P$-bounded with relative bound less than one. For $|\alpha|+|\beta|\leq 2$, $D_\rx^\alpha D_\ry^\beta$ is $\Delta$-bounded. The above properties imply that it is also $(P+W)$-bounded, yielding \eqref{eq:borne-elliptique}. \cqfd

For $\epsilon\geq 0$, let $\Omega _\epsilon :=\{x\in\Omega _0;\, d(x;\R^{3m}\setminus\Omega _0)>\epsilon\}$ and, for $r\in\N$, denote the $\rL^2(\Omega _\epsilon ;\cW _r)$-norm of $v$ by $N_{\epsilon ; r}(v)$. Note that $\chi _0\varphi\in\cC_c^\infty(\R^{3m}; \cW_2)$ with support in $\Omega _0'$, by Lemma~\ref{l:smooth}. 
\begin{lemma}\label{l:anal-control}
There exists $B>0$ such that, for all $\epsilon>0$, $j\in\N$, $r\in\{0;1;2\}$, and $\alpha\in\N^{3m}$,
\begin{equation}\label{eq:anal-control}
r+|\alpha |<2+j\ \impl \epsilon ^{r+|\alpha |}N_{j\epsilon ; r}(D_x^\alpha\varphi)\ \leq \ B^{r+|\alpha |+1}\period
\end{equation}
\end{lemma}
\Pf Let $D_0$ be the diameter of $\Omega _0$. Since $P\in{\rm Diff}_2(\Omega)$, there exists $C_a>0$ such that, for all $\alpha\in\N^3$, $0\leq \epsilon_1\leq D_0$, 
\begin{equation}
\epsilon_1^{|\alpha|}\sum_{|\beta |\leq 2}\, \sup _{x\in\Omega _{\epsilon_1}}\, \|\partial _x^\alpha a_\beta \|_{\cB_{2-|\beta |}}\ \leq \ C_a^{|\alpha |+1}\cdot (|\alpha |!)\period\label{eq:coeff-anal}
\end{equation}
by \eqref{eq:anal-2}. Using that $(P+W)\varphi =0$, it suffices to follow the proof of (7.5.9) in \cite{h1}, p. 179-180, to get the desired result. \cqfd

To complete the proof of Theorem~\ref{th:ellipti-regu}, it suffices, 
thanks to \eqref{eq:anal-control}, to follow the argument on page 180 in \cite{h1} and use \eqref{eq:anal-2}. \cqfd

\begin{remark}\label{r:sans-pseudo}
We point out that one can prove Theorem~\ref{th:ellipti-regu} without 
pseudodifferential technics. We used them in the proof of Lemma~\ref{l:smooth}. In fact, it turns out that the arguments used in {\rm\cite{fhhs2}} {\rm (}see the proof of Lemma 3.1{\rm )} and in {\rm\cite{hs}} {\rm (}see Section 3.3{\rm )} can be adapted to get a different proof of those lemma. 
\end{remark}
%

\section{Twisted differential calculus.}
\label{s:twisted-diff}
\setcounter{equation}{0}

In this section, we recall the notion of ``twist'' used in \cite{hu,kmsw,ms} to ``desingularize'' a class of Hamiltonians with Coulomb-like singularities.

\subsection{Framework.}
\label{ss:cadre}

We need some more notation. Denoting by $\SSS^2$ the unit euclidian sphere of $\R^3$, let $\tilde{v}_0 : \SSS^2\dans\C$ and $\eta : ]0; +\infty[\dans\R^+$ be two functions such that, for some $\eta_0\in ]0; 1[$ and $B_0>0$, $(1/B_0)\leq |\tilde{v}_0|\leq B_0$ and 
\begin{equation}\label{eq:cond-eta}
 \forall t>0\comma\hspace{.4cm}\sup_{t(1-\eta _0)\, \leq\, s\, \leq\, t(1+\eta_0)}\eta (s)\ \leq \ B_0\cdot \eta (t)\period
\end{equation}
We assume further that the multiplication operator by the function $\eta (|\cdot|) : \R^3\ni z\donne \eta (|z|)\in\R^+$ belongs to $\cL (\cW^{1,2}(\R^3); \rL^2(\R^3))$. The multiplication operator 
by the function $v_0 : \R^3\setminus\{0\}\dans\C$, defined by $v_0(z)=\eta (|z|)\tilde{v}_0(z/|z|)$, also belongs to $\cL (W^{1,2}(\R^3); \rL^2(\R^3))$. 
\begin{definition}\label{def:classe-V}
 Let $\cV$ be the class of functions $v: \R^3\dans\C$ such that 
$v\in\cC^\omega(\R^3\setminus\{0\})$ and there 
exists $C>0$ such that
\begin{equation}\label{eq:cond-v}
 \forall \alpha\in\N^3\comma\ \forall z\in\R^3\setminus\{0\}\comma\ |z|^{|\alpha|}\cdot\bigl|D^\alpha v(z)\bigr|\ \leq \ C^{1+|\alpha|}\cdot \bigl(|\alpha|!\bigr)\cdot |v_0(z)|\period
\end{equation}
\end{definition}
We observe that, if $v\in\cV$, then the multiplication operator by $v$ belongs to the space $\cL (\cW^{1,2}(\R^3); \rL^2(\R^3))$. 

We consider a system of $m+p$ electrons, interacting to one another and moving under the influence of $L$ fixed nuclei. Let $\cI _e:=\{1; \cdots ; m\}$ be the set of indices for the ``external'' electronic variables $\rx\in\R^{3m}$, let $\cI _i:=\{1; \cdots ; p\}$ be the one of ``internal'' electronic variables $\ry\in\R^{3p}$, and let $\cI _n:=\{1; \cdots ; L\}$ be the one of the nuclear variables $R_\ell\in\R^3$. It is convenient to introduce the disjoint union of those sets, namely 
\[\cI _e\, \sqcup\, \cI _i\, \sqcup\, \cI _n\ :=\ \bigl\{(a; j)\, ;\ a\in\{e; i; n\}\comma\ j\in\cI_a\bigr\}\period\]
For $c=(a; j)\in \cI _e\sqcup\cI _i\sqcup\cI _n$, let $\pi_1(c)=a$ and $\pi_2(c)=j$. We define $X_c=\rx_j$ if $\pi_1(c)=e$, $X_c=\ry_j$ if $\pi_1(c)=i$, and $X_c=R_j$ if $\pi_1(c)=n$. 
Let 
\[\cP \ :=\ \bigl\{\{c; c'\}\, ;\ (c; c')\in \bigl(\cI _e\, \sqcup\, \cI _i\, \sqcup\, \cI _n\bigr)^2\ \mbox{and}\ c\neq c'\bigr\}\]
be the set of all possible particle pairings. Let $\uOmega$ be an open subset of $\R^{3m}$. We consider potentials $\uV$ where 
\begin{equation}\label{eq:def-V}
\forall (\rx ; \ry)\in \uOmega\times\R^{3p}\comma\hspace{.4cm}\uV (\rx ; \ry)\ =\ \sum _{\{c;\, c'\}\in\cP}\, V_{\{c;\, c'\}}(X_c\, -\, X_{c'})\comma
\end{equation}
with $V_{\{c;\, c'\}}\in\cV$ (cf. Definition~\ref{def:classe-V}), for all $\{c; c'\}\in\cP$. \\
Recall that $\cW_n=W^{n,2}(\R^{3p})$, for $n\in\N$, and $\cB_1=\cL (\cW_1; \cW_0)$. We note that, if $v\in\cV$ and 
$\{c; c'\}\in\cP$, then, for any fixed $\rx$, the multiplication by the almost everywhere defined function $v(X_c-X_{c'})$ of $\ry$ belongs to $\cB_1$, thanks to \eqref{eq:cond-v} and the assumption on $v_0$. 

Recall that we introduced in Definition~\ref{def:ope-diff} a class of differential operators. The class of Hamiltonians, that we want to treat, is defined in the following:
\begin{definition}\label{def:hamiltonians}
 For $(m; p)\in (\N^\ast)^2$ and $\uOmega$ an open subset of $\R^{3m}$. We denote by ${\rm Ham}_2(\uOmega)$ the set of all the Hamiltonians $\uH$ on $\uOmega\times\R^{3p}$ such that $\uH=\uP+\uV$, where $\uP\in {\rm Diff}_2(\Omega)$ {\rm (}cf. Definition~\ref{def:ope-diff}{\rm )} and $\uV$ satisfies \eqref{eq:def-V}.
\end{definition}
%

\subsection{Heuristic.}
\label{ss:heuristic}

In this subsection, we want to motivate the notion of ``twist'', that is going to play a crucial rôle in the proof of the main results. To explain the basic idea, let us consider the following situation. 

Take the Hamiltonian $H$ in \eqref{eq:hamiltonien} with $N=2$. We are not able to apply Theorem~\ref{th:ellipti-regu} with $\rx =x_1$ and $\ry =x_2$ to the equation $H\psi =E\psi$, since the map $x_1\donne (x_2\donne |x_1-x_2|^{-1})\in\cB _1=\cL (W^{1,2}(\R^3); \rL^2(\R^3))$ is admittedly well defined by \eqref{eq:hardy} but is not analytic. The problem lies in the dependency of the Coulomb singularity on $x_1$. \\
Assume that, for some open subset $\Omega$ of $\R^3$ and some $x^0\in\R^3$, we have a smooth map $f : \Omega\times\R^3\to \R^3$ such that, for all $x_1\in\Omega$, $f(x_1; \cdot)$ is a smooth diffeomorphism of $\R^3$ satisfying $f(x_1; x^0)=x_1$. For such $x_1$, let $U_{x_1}$ be the unitary implementation of $f(x_1; \cdot)$ on $\rL^2(\R^3)$ (the ``twist''). For any function $\varphi\in W^{1,2}(\R^3)$ and 
for $x_1\in\Omega$, 
\[U_{x_1}\bigl(|x_1\, -\, \cdot|^{-1}\varphi \bigr)(y)\ =\ \bigl|f(x_1; x^0)\, -\, f(x_1; y)\bigr|^{-1}\cdot \bigl(U_{x_1}\varphi\bigr)(y)\period\]
Now, since $f(x_1; \cdot)$ is a diffeomorphism, the singularity occurs when $y=x^0$. Therefore, it does not depends on $x_1$ anymore. We can write 
\[\bigl(U_{x_1}|x_1\, -\, \cdot|^{-1}U_{x_1}^{-1}\varphi \bigr)(y)\ =\ |x^0-y|\cdot \bigl|f(x_1; x^0)\, -\, f(x_1; y)\bigr|^{-1}\cdot |x^0-y|^{-1}\cdot \varphi (y)\period\]
Assuming further that the map 
\begin{equation}\label{eq:bornitude}
 \Omega\ni x_1\, \donne\, |x^0-\cdot|\cdot \bigl|f(x_1; x^0)\, -\, f(x_1; \cdot)\bigr|^{-1}\in\rL^\infty (\R^3)
\end{equation}
is smooth, we see that $\Omega\ni x_1\donne U_{x_1}|x_1-\cdot|^{-1}U_{x_1}^{-1}\in\cB_1=\cL(W^{1,2}(\R^3); \rL^2(\R^3))$ is also smooth, thanks to \eqref{eq:hardy}.  \\
What can we do with the other singular terms in $H$? Since the terms $|x_1-R_\ell|^{-1}$ do not depend on $x_2$, the previous conjugation has no effect on them. But, if $x_1$ stays far away from all $R_\ell$, they are not singular. We are left with the terms $|x_2-R_\ell|^{-1}$, which are a constant function of $x_1$ with values in $\cB_1$ by \eqref{eq:hardy}. Since we want to perform the above conjugation to treat the term $|x_1-x_2|^{-1}$, we need to verify that this conjugation does not distroy the nice property of these terms. If we require that, for all $\ell$ and all $x_1\in\Omega$, $f(x_1; R_\ell)=R_\ell$, and the smoothness of the map \eqref{eq:bornitude} with $x^0$ replaced by $R_\ell$, the above computation shows that the map $\Omega\ni x_1\donne U_{x_1}|x_2-R_\ell|^{-1}U_{x_1}^{-1}\in\cB_1$ is also smooth. \\
With some more work, we can prove real analyticity of the previous terms. 
We just have seen that the ``twisted'' potential $U_{x_1}VU_{x_1}^{-1}$ is a real analytic, $\cB_1$-valued function. What about the ``twisted'' kinetic energy? A simple computation shows that the operators $\nabla _{x_1}$ and $\nabla _{x_2}$ are transformed by the twist into differential operators in the $(x_1; x_2)$ variables with explicit coefficients in terms of the derivatives of $f$ (cf. \eqref{eq:twisted-1-gradiant-x} and \eqref{eq:twisted-1-gradiant-y}). As we shall see below, the ellipticity of $H$ is preserved by the ``twist'' and the coefficients of the ``twisted'' operator have some property of real analyticity. \\
Coming back to the original problem, we do not apply Theorem~\ref{th:ellipti-regu} with $\rx =x_1$ and $\ry =x_2$ to the equation $H\psi =E\psi$ but to the ``twisted'' equation $(U_{x_1}HU_{x_1}^{-1})U_{x_1}\psi =EU_{x_1}\psi$. Therefore, we obtain some analytic regularity of the ``twisted'' bound state, namely $U_{x_1}\psi$, not on $\psi$ itself. But, if we are interested in the 
regularity of $\rho _1$, we are done. Indeed, since we can write $\rho _1(x_1)=\|\psi (x_1; \cdot )\|^2$ ($\|\cdot\|$ being the $\rL^2$-norm on $\R^3$) and since $U_{x_1}$ is unitary on $\rL^2(\R^3)$, $\rho _1(x_1)=\|U_{x_1}\psi (x_1; \cdot )\|^2$ and $\rho _1$ is real analytic on $\Omega$ by composition. \\
We just have sketched, in the $N=2$ case, the proof of the smoothness of $\rho _1$, which was performed in \cite{j}, provided there exists a map $f$ satisfying the above requirements. Actually it does (see Subsection~\ref{ss:twist} below).

\subsection{Hunziker's twist.}
\label{ss:twist}

With the above strategy in mind, we recall Hunziker's twist (cf. \cite{hu,kmsw}) and provide basic properties of it.

We take a real-valued cut-off function $\tau\in\cC_c^\infty(\R^3)$ such that $\tau (x)=0$ if $|x|\geq 1$, and $\tau (0)=1$. Let $\cU$ be the open subset of $\R^{3m}$ defined by 
\begin{align}
 \cU\ =\ & \bigl\{\rx\in\R^{3m}\, ;\ \forall (c; c')\in\cP;\ \pi_1(c)\, =\, e\ \mbox{and}\ \pi_1(c')\in\{e; n\}\comma\hspace{.4cm}X_c\neq X_{c'}\bigr\}\nonumber\\
 =\ &\bigl\{\rx\in\R^{3m}\, ;\ \forall (j; j')\in\ncg 1; m\ncd^2\comma\ \forall \ell\in\ncg1; L\ncd\comma\ \rx_j\neq R_\ell\label{eq:def-U-ronde}\\
 &\hspace{7.4cm}\mbox{and}\ j\neq j'\, \impl\, \rx_j\neq\rx_{j'}\bigr\}\period\nonumber
\end{align}
Let $\rx^0=(\rx_1^0 ; \cdots ; \rx_m^0)\in\cU$. We can find some $r_0>0$ such that, for $(j; j')\in \ncg 1; m\ncd^2$, for $\ell\in\ncg1; L\ncd$,  
\begin{equation}\label{eq:écartement}
 |\rx_j^0-R_\ell|\ \geq\ r_0\hspace{.4cm}\mbox{and}\hspace{.4cm}j\neq j'\, \impl\, |\rx_j^0\, -\, \rx_{j'}^0|\ \geq\ r_0\period
\end{equation}
Let $f : \R^{3m}\times\R^3\dans\R^3$ be the smooth function defined by 
\begin{equation}\label{eq:def-f}
 f(\rx ; z)\ =\ z\, +\, \sum _{j=1}^m\, \tau\bigl(r_0^{-1}(z\, -\rx_j^0)\bigr)\cdot (\rx _j\, -\, \rx _j^0)\period
\end{equation}
Notice that, by \eqref{eq:écartement}, we have, for all $\rx\in\R^{3m}$, 
\begin{equation}\label{eq:transformation}
\forall j\in\ncg1; m\ncd\comma\ \ell\in\ncg 1; L\ncd\comma\hspace{.4cm}f(\rx ; \rx_j^0)\ =\ \rx _j
\hspace{.4cm}\mbox{and}\hspace{.4cm}f(\rx ; R_\ell)\ =\ R_\ell\period
\end{equation}
For all $(\rx ; z)\in\R^{3m}\times\R^3$, the total differential of $f$ w.r.t. $\rx$ is given by 
\begin{equation}\label{eq:dxf}
 (d_\rx f)(\rx; z)\ =\ \sum _{j=1}^m\, \tau\bigl(r_0^{-1}(z\, -\rx_j^0)\bigr)\ d\rx _j
\end{equation}
where $d\rx _j : \R^{3m}\ni\rx '\donne\rx _j'\in\R^3$, for $1\leq j\leq m$. In particular, 
the map $\R^{3m}\ni\rx\donne (d_\rx f)(\rx; \cdot)\in\cC_b^\infty(\R^3; \cL (\R^{3m}; \R^3))$ is well-defined and constant.
Using \eqref{eq:dxf} and setting 
\[C(\tau )\ =\ \sup_{z\in\R^3}\, \bigl\|d\tau (z)\bigr\|_{\cL (\R^3)}\ <\ +\infty\comma\]
we can write, for all $j\in\ncg1; m\ncd$, for all $(z; z')\in (\R^3)^2$, and for all $\rx\in\R^{3m}$, 
\begin{align}
 \bigl\|(d_{\rx _j} f)(\rx ; z)\, -\, (d_{\rx _j} f)(\rx ; z')\bigr\|_{\cL (\R^{3m}; \R^3)}\ =&\ \Bigl|\tau \bigl(r_0^{-1}(z\, -\rx_j^0)\bigr)\, -\, \tau \bigl(r_0^{-1}(z'\, -\rx_j^0)\bigr)\Bigr|\nonumber\\
 \leq& \ 
 r_0^{-1}C(\tau )\, |z\, -\, z'|\period\label{eq:lipschitz-1}
\end{align}
For all $(\rx ; z)\in\R^{3m}\times\R^3$, the total differential of $f$ w.r.t. $z$ is given by, for $z^\epsilon\in\R^3$,  
\begin{equation}\label{eq:dzf}
 (d_zf)(\rx; z)\cdot z^\epsilon\ =\ z^\epsilon\, +\, r_0^{-1}\sum _{j=1}^m\, \bigl((d\tau) \bigl(r_0^{-1}(z\, -\rx_j^0)\bigr)\cdot z^\epsilon\bigr)(\rx _j\, -\, \rx _j^0)\period
\end{equation}
For $\delta>0$, set $\Omega (\delta ):=B(\rx_1^0; \delta [\times\cdots\times B(\rx_m^0; \delta [$. Recall that the class $\cV$ (cf. Definition~\ref{def:classe-V}) depends on a parameter $\eta_0$ (see \eqref{eq:cond-eta}). 
By \eqref{eq:dzf}, we can find some $\delta _0\in ]0; r_0/2]$, depending on $C(\tau )$ and $r_0$, such that
\begin{equation}\label{eq:dzf-1}
 \forall\rx\in\Omega (\delta _0)\comma\hspace{.4cm} \sup _{z\in\R^3}\, \bigl\|(d_zf)(\rx; z)\, -\, I_3\bigr\|_{\cL (\R^3)}\ \leq\ \min \bigl(\eta_0;\, 1-\eta_0\bigr)\ <\ 1\period
\end{equation}
We note that the requirement $\delta _0\in ]0; r_0/2]$ yields $\Omega (\delta _0)\subset\cU$, by \eqref{eq:écartement} and the triangle inequality. Thanks to \eqref{eq:dzf-1}, for each $\rx\in\Omega (\delta _0)$, $f(\rx ; \cdot)$ is a $\cC^\infty$-diffeomorphism of $\R^3$. We denote by $f^{\langle -1\rangle}(\rx; \cdot)$ its inverse.
Furthermore, for all $\rx\in\Omega (\delta _0)$, for all $(z; z')\in (\R^3)^2$,  
\begin{equation}\label{eq:lipschitz-2}
 (1\, -\, \eta_0)\, |z\, -\, z'|\ \leq \bigl|f(\rx ; z)\, -\, f(\rx ; z')\bigr|\ \leq \ 
 (1\, +\, \eta_0)\, |z\, -\, z'|\period
\end{equation}
By a direct computation using \eqref{eq:dxf} and \eqref{eq:dzf}, one can express the derivative $d_zf^{\langle -1\rangle}$ and $d_\rx f^{\langle -1\rangle}$ in terms of $d_zf$ (see Appendix~\ref{app:computations} for details). In particular, one can check that $f\in\cC_b^\infty (\Omega (\delta _0)\times\R^3; \R^3)$ and $f^{\langle -1\rangle}\in\cC_b^\infty (\Omega (\delta _0)\times\R^3; \R^3)$. \\
Consider the map $F : \Omega (\delta _0)\times\R^{3p}\dans\R^{3p}$, defined by $F(\rx; \ry)=(f(\rx; \ry_1); \cdots ; f(\rx; \ry_p))$. For each $\rx\in\Omega (\delta _0)$, $F(\rx; \cdot)$ is a $\cC^\infty$-diffeomorphism of $\R^{3p}$ and we denote its inverse map by $F^{\langle -1\rangle}(\rx; \cdot)$. Explicit expressions for $d_\rx F$, $d_\ry F$, $d_\rx F^{\langle -1\rangle}$, and $d_\ry F^{\langle -1\rangle}$, are provided in Appendix~\ref{app:computations}. In particular, we observe that $F\in\cC_b^\infty (\Omega (\delta _0)\times\R^{3p}; \R^{3p})$ and 
$F^{\langle -1\rangle}\in\cC_b^\infty (\Omega (\delta _0)\times\R^{3p}; \R^{3p})$. \\
We introduce the map $U : \Omega (\delta _0)\ni \rx\donne U_\rx\in\cB_0$ (with $\cB_0=\cL (\cW _0)$ and $\cW _0=\rL^2(\R^{3p})$), with values in the set of unitary operators on $\cW_0$, defined by, for $\theta\in\cW_0$, 
\begin{equation}\label{eq:def-U}
 \bigl(U_\rx\theta\bigr)(\ry)\ =\ \bigl|{\rm Det}(d_\ry F)(\rx; \ry)\bigr|^{1/2}\cdot \theta\bigl(F(\rx; \ry)\bigr)\period
\end{equation}
For $\rx\in\Omega (\delta _0)$, the map $U_\rx$ is the unitary implementation of $F(\rx ; \cdot)$ on $\rL^2(\R^{3p})$ and 
\begin{equation}\label{eq:def-U-bis}
 \bigl(U_\rx^{-1}\theta\bigr)(\ry)\ =\ \Bigl|{\rm Det}\bigl(d_\ry F^{\langle -1\rangle}\bigr)(\rx; \ry)\Bigr|^{1/2}\cdot \theta\bigl(F^{\langle -1\rangle}(\rx; \ry)\bigr)\period
\end{equation}
Denoting by $A^T$ the transposed of a linear map $A$ and by ${\rm Det}A$ its determinant, we see by a direct computation that, as differential operators on $\Omega (\delta _0)\times\R^{3p}$, 
\begin{align}
U_{\rx}^{-1}\, D_\rx\, U_{\rx}\ =&\ D_\rx\, +\, J_1(F)\, D_\ry\, +\, J_2(F)\comma\label{eq:twisted-1-gradiant-x}\\
U_{\rx}\, D_\rx\, U_{\rx}^{-1}\ =&\ D_\rx\, +\, J_1\bigl(F^{\langle -1\rangle}\bigr)\, D_\ry\, +\, J_2\bigl(F^{\langle -1\rangle}\bigr)\comma\label{eq:twisted-2-gradiant-x}\\
U_{\rx}^{-1}\, D_{\ry '}\, U_{\rx}\ =&\ J_3(F)\, D_\ry\, +\, J_4(F)\comma\label{eq:twisted-1-gradiant-y}\\
U_{\rx}\, D_{\ry '}\, U_{\rx}^{-1}\ =&\ J_3\bigl(F^{\langle -1\rangle}\bigr)\, D_\ry\, +\, J_4\bigl(F^{\langle -1\rangle}\bigr)\comma\label{eq:twisted-2-gradiant-y}
\end{align}
where, for $\rx\in\Omega (\delta _0)$ and $\ry\in\R^{3p}$, with $G=F$ and $G^{\langle -1\rangle}=F^{\langle -1\rangle}$ or with $G=F^{\langle -1\rangle}$ and $G^{\langle -1\rangle}=F$, 
\begin{align*}
 J_1(G)(\rx ; \ry)\ :=&\ \Bigl(\bigl(d_\rx G\bigr)(\rx ; \ry')\Bigr)^{\mathrm{T}}\bigl(\rx ; \, \ry'=G^{\langle -1\rangle}(\rx ; \ry)\bigr)\comma \\
J_2(G)(\rx ; \ry)\ :=&\ -i\bigl|{\rm Det}\, \bigl(d_\ry G^{\langle -1\rangle}\bigr)(\rx ; \ry)\bigr|^{1/2}\, D_\rx\Bigl(\bigl|{\rm Det}\, \bigl(d_{\ry'}G\bigr)(\rx ; \ry ')\bigr|^{1/2}\Bigr)\Bigr|_{\ry'=G^{\langle -1\rangle}(\rx ; \ry)}\comma\\
J_3(G)(\rx ; \ry)\ :=&\ \Bigl(\bigl(d_{\ry'}G\bigr)(\rx ; \ry ')\Bigr)^{\mathrm{T}}\bigl(\rx; \, \ry'=G^{\langle -1\rangle}(\rx ; \ry)\bigr)\comma \\
J_4(G)(\rx ; \ry)\ :=&\ -i\bigl|{\rm Det}\, \bigl(d_\ry G^{\langle -1\rangle}\bigr)(\rx ; \ry)\bigr|^{1/2}\, D_{\ry '}\Bigl(\bigl|{\rm Det}\, \bigl(d_{\ry'}G\bigr)(\rx ; \ry ')\bigr|^{1/2}\Bigr)\Bigr|_{\ry'=G^{\langle -1\rangle}(\rx ; \ry)}\period 
\end{align*}
Thanks to the computation written in Appendix~\ref{app:computations}, we note that, on one hand, for $G\in\{F; F^{\langle -1\rangle}\}$,
\begin{align}
 J_1(G)\ \in\ \cC_b^\infty \Bigl(\Omega (\delta _0)\times\R^{3p};\, \cL\bigl(\R^{3p}; \R^{3m}\bigr)\Bigr)\comma&\hspace{.4cm}iJ_2(G)\ \in\ \cC_b^\infty \bigl(\Omega (\delta _0)\times\R^{3p};\, \R^{3m}\bigr)\comma\label{eq:prop-J-1}\\
 J_3(G)\ \in\ \cC_b^\infty \Bigl(\Omega (\delta _0)\times\R^{3p};\, \cL\bigl(\R^{3p}\bigr)\Bigr)\comma&\hspace{.4cm}iJ_4(G)\ \in\ \cC_b^\infty \bigl(\Omega (\delta _0)\times\R^{3p};\, \R^{3p}\bigr)\comma\label{eq:prop-J-2}
\end{align}
and, on the other hand, if we denote by $J$ the multiplication operator by any component of $J_k(G)$, $k\in\{1; 2; 3; 4\}$ and $G\in\{F; F^{\langle -1\rangle}\}$, 
\begin{equation}\label{eq:prop-anal-J}
 J\ \in\ \cC^\omega\bigl(\Omega (\delta _0);\, \cB_0\bigr)\period
\end{equation}
Furthermore (cf. Appendix~\ref{app:computations}), there exists $C>0$ such that, for $G\in\{F; F^{\langle -1\rangle}\}$, for all $(\rx ; \ry)\in\Omega (\delta _0)\times\R^{3p}$, $J_3(G)(\rx ; \ry )$ is invertible and 
\begin{equation}\label{eq:inv-J3}
 C^{-1}\ \geq \ \bigl\|J_3(G)(\rx ; \ry)\bigr\|_{\cL (\R^{3p})}\, +\, \bigl\|\bigl(J_3(G)(\rx ; \ry)\bigr)^{-1}\bigr\|_{\cL (\R^{3p})}\ \geq\ C\period 
\end{equation}
Now, using \eqref{eq:twisted-1-gradiant-y}, 
\[D_{\ry '}U_\rx\ =\ U_\rx\bigl(J_3(F)D_\ry - J_4(F)\bigr)\]
thus, for all $\rx\in\Omega (\delta _0)$, we see, by induction and by \eqref{eq:prop-J-2}, that $U_\rx$ leaves $\cW_n$ invariant, for all $n\in\N$.
Similarly, we can show, using \eqref{eq:twisted-2-gradiant-y} and \eqref{eq:prop-J-2}, that, for all $\rx\in\Omega (\delta _0)$, $U_\rx^{-1}$ also leaves $\cW_n$ invariant, for all $n\in\N$. By the same arguments, we can check, starting from \eqref{eq:twisted-1-gradiant-x} and \eqref{eq:twisted-1-gradiant-y}, that 
$U$ preserves the space $W^{2,2}(\Omega (\delta _0)\times \R^{3p})$.

\subsection{Potential regularization by a twist.}
\label{ss:twisted-potential}

We first consider the action by conjugation of the twist $U$ on the potential part $\uV$ of Hamiltonians $\uH$ in the class ${\rm Ham}_2(\uOmega)$. 
\begin{proposition}\label{p:pot-regu-anal} 
Take $\rx^0\in\cU\subset\R^{3m}$. We choose $r_0>0$ and $\delta _0\in ]0; r_0/2]$ such that \eqref{eq:écartement} and \eqref{eq:dzf-1} hold true. Recall that the neighbourhood $\Omega (\delta_0):=B(\rx_1^0; \delta _0[\times\cdots\times B(\rx_m^0; \delta _0[$ of $\rx^0$ is included in $\cU$. \\
Then, for any $v\in\cV$ and any $\{c; c'\}\in\cP$, the map $\Omega (\delta_0)\ni\rx\donne U_\rx v(X_c-X_{c'})U_\rx ^{-1}\in\cB_1$ is well-defined and real analytic. 
\end{proposition}
\begin{remark}\label{r:coulomb-twisté}
 We point out that the Coulomb potential $|\cdot|^{-1}$ does belong to $\cV$. 
 Indeed, \eqref{eq:cond-v}, for $v=v_0=|\cdot|^{-1}$, was proved in {\rm \cite{j}} and the multiplication by this $v_0$ belongs to $\cL (\cW^{1,2}(\R^3); \rL^2(\R^3))$ by Hardy's inequality \eqref{eq:hardy}. 
\end{remark}
\begin{remark}\label{r:pot-twisté}
 A proof of Proposition~\ref{p:pot-regu-anal} is essentially available in {\rm \cite{hu,j}}. However, the proof in {\rm \cite{hu}} is restricted to Coulomb interactions. We provide below a more transparent version of the proof in {\rm \cite{j}}. 
\end{remark}
{\bf Proof of Proposition~\ref{p:pot-regu-anal}:} Let $\rx\in\Omega (\delta_0)$. Recall that, for $k\in\{0; 1\}$, $U_\rx$ and $U_\rx^{-1}$ leave $\cW_k$ invariant. Recall that, for all $\{c; c'\}\in\cP$, the multiplication by $v(X_c-X_{c'})$ belongs to $\cB_1$. Therefore the map $w: \Omega (\delta_0)\ni\rx\donne U_\rx v(X_c-X_{c'})U_\rx ^{-1}\in\cB_1$ is well-defined. We treat all possible cases for $\{c; c'\}\in\cP$. \\
If $\pi _1(c)=\pi _1(c')=n$, then $w : \Omega (\delta_0)\ni\rx\donne v(X_c-X_{c'})$ is a constant. \\
Assume that $\pi _1(c)=e$ and $\pi _1(c')\in\{e; n\}$ or that $\pi _1(c')=e$ and $\pi _1(c)\in\{e; n\}$. Since $\Omega (\delta_0)\subset\cU$, $X_c-X_{c'}$ does not vanish on $\Omega (\delta_0)$. Thus 
$w : \Omega (\delta _0)\ni\rx\donne v(X_c-X_{c'})$ is real analytic, by composition. \\
Assume that $\pi _1(c)=e$ and $\pi _1(c')=i$. Let $j=\pi _2(c)$ and $k=\pi _2(c')$. Then, for $\rx\in\Omega (\delta_0)$, 
$U_\rx v(X_c-X_{c'})U_\rx ^{-1}$ is the multiplication operator on $\cW _1$ by the function 
\[\ry\ \donne\ v\bigl(\rx _j\, -\, f(\rx; \ry _k)\bigr)\ =\ v\bigl(f(\rx; \rx_j^0)\, -\, f(\rx; \ry _k)\bigr)\comma\]
by \eqref{eq:transformation}. Let $j'\in\cI _e$. Using the fact that $d_\rx f$ does not depend on $\rx$ (see \eqref{eq:dxf}), we have, for $\alpha\in\N^{3}$ and $\ry_k\neq\rx_j^0$, 
\begin{align}\label{eq:dérive-usuelle}
&\partial _{\rx _{j'}}^\alpha\Bigl(v\bigl(f(\rx; \rx_j^0)\, -\, f(\rx; \ry _k)\bigr)\Bigr)\\
=\ &(\partial^\alpha v)\bigl(f(\rx; \rx_j^0)\, -\, f(\rx; \ry _k)\bigr)\, \cdot\, \Bigl(\tau \bigl(r_0^{-1}(\rx_j^0\, -\rx _{j'}^0)\bigr)\, -\, \tau \bigl(r_0^{-1}(\ry _k\, -\rx _{j'}^0)\bigr)\Bigr)^{|\alpha|}\period
\nonumber
\end{align}
Using successively \eqref{eq:lipschitz-1}, \eqref{eq:lipschitz-2}, and \eqref{eq:cond-v}, we obtain 
\begin{align*}
&\Bigl|\partial _{\rx _{j'}}^\alpha\Bigl(v\bigl(f(\rx; \rx_j^0)\, -\, f(\rx; \ry _k)\bigr)\Bigr)\Bigr|\\
\leq\ &\bigl|(\partial^\alpha v)\bigl(f(\rx; \rx_j^0)\, -\, f(\rx; \ry _k)\bigr)\bigr|\, \cdot\, \Bigl(r_0^{-1}\, C(\tau)\, |\rx_j^0\, -\, \ry_k|\Bigr)^{|\alpha|}\comma\\
\leq\ &\bigl|(\partial^\alpha v)\bigl(f(\rx; \rx_j^0)\, -\, f(\rx; \ry _k)\bigr)\bigr|\, \cdot\, \Bigl(\bigl(r_0(1-\eta _0)\bigr)^{-1}\, C(\tau)\, \bigl|f(\rx; \rx_j^0)\, -\, f(\rx; \ry _k)\bigr|\Bigr)^{|\alpha|}\comma\\
\leq\ &C_1^{1+|\alpha|}\cdot \bigl(|\alpha|!\bigr)\cdot \Bigl|v_0\bigl(f(\rx; \rx_j^0)\, -\, f(\rx; \ry _k)\bigr)\Bigr|\comma
\end{align*}
for some $C_1>0$. Now, by the structure of $v_0$ and by \eqref{eq:lipschitz-2} and \eqref{eq:cond-eta}, 
we derive that, for some positive constants $C_2$ and $C_3$, 
\begin{align}
\Bigl|\partial _{\rx _{j'}}^\alpha\Bigl(v\bigl(f(\rx; \rx_j^0)\, -\, f(\rx; \ry _k)\bigr)\Bigr)\Bigr|
&\leq\ C_2^{1+|\alpha|}\cdot \bigl(|\alpha|!\bigr)\cdot \eta\bigl(f(\rx; \rx_j^0)\, -\, f(\rx; \ry _k)\bigr)\comma\nonumber\\
&\leq\ C_3^{1+|\alpha|}\cdot \bigl(|\alpha|!\bigr)\cdot\eta\bigl(|\rx_j^0\, -\, \ry_k|\bigr)\period\label{eq:anal-bound}
\end{align}
Now, one can check, by induction on $|\alpha|$ and by the fact that the multiplication by $\eta (|\cdot|)$ belongs to $\cB_1$, that, on $\Omega (\delta_0)$, the $\partial _{\rx _{j'}}^\alpha$ partial derivative of $w$ exists and maps each $\rx\in\Omega (\delta_0)$ to the multiplication operator on $\cW_1$ by the almost everywhere defined function \eqref{eq:dérive-usuelle} of $\ry$. Thus, $w\in\cC^\infty(\Omega (\delta_0); \cB_1)$. Using \eqref{eq:anal-bound} again, we obtain that 
\begin{equation}\label{eq:bound-derivative-v}
 \sup _{\rx\in\Omega (\delta_0)}\, 
\bigl\|(D_{\rx_{j'}}^\alpha w)(\rx)\bigr\|_{\cB_1}\ \leq \ C_4^{1+|\alpha|}\cdot \bigl(|\alpha|!\bigr)\comma
\end{equation}
for some $C_4>0$. By \eqref{eq:anal-2}, $w$ is real analytic w.r.t. $\rx_{j'}$. Since this holds true for all $j'$, $w$ is real analytic. \\
Similarly, we can treat the case $\pi _1(c)=i$ and $\pi _1(c')=e$.\\
Assume that $\pi _1(c)=n$ and $\pi _1(c')=i$. Let $\ell=\pi _2(c)$ and $k=\pi _2(c')$. For $\rx\in\Omega (\delta _0)
$, $w(\rx)$ is the multiplication operator in $\cB_1$ by the function 
\begin{equation}\label{eq:v-R-l}
 \ry\ \donne\ v\bigl(R_\ell\, -\, f(\rx; \ry _k)\bigr)\ =\ v\bigl(f(\rx; R_\ell)\, -\, f(\rx; \ry _k)\bigr)\comma
\end{equation}
by \eqref{eq:transformation}. We redo the computation \eqref{eq:dérive-usuelle} - \eqref{eq:anal-bound}, with $\rx_j^0$ replaced by $R_\ell$. Furthermore, we also can check that the $\partial _{\rx _{j'}}^\alpha$ partial derivative of $w$ exists and maps each $\rx\in\Omega (\delta_0)$ to the multiplication operator on $\cW_1$ by the almost everywhere defined $\partial _{\rx _{j'}}^\alpha$ partial derivative of \eqref{eq:v-R-l}.
Using the new estimate \eqref{eq:anal-bound}, we derive a bound like 
\eqref{eq:bound-derivative-v}. As above, we conclude that $w$ is real analytic. \\
Again, we have a similar treatment of the case $\pi _1(c)=i$ and $\pi _1(c')=n$.\\
We are left with the case $\pi _1(c)=i$ and $\pi _1(c')=i$. We set $k=\pi _2(c)$ and $k'=\pi _2(c')$. By assumption, $k\neq k'$. Using \eqref{eq:transformation} again, we see that, for $\rx\in\Omega (\delta _0)
$, $w(\rx)$ is the multiplication operator in $\cB_1$ by the function 
\begin{equation}\label{eq:v-y-y}
 \ry\ \donne\ v\bigl(f(\rx; \ry_k)\, -\, f(\rx; \ry _{k'})\bigr)\period
\end{equation}
Once again, we perform the computation \eqref{eq:dérive-usuelle} - \eqref{eq:anal-bound} and follow the above arguments, using this time that the multiplication by $\ry\donne \eta (|\ry_k-\ry _{k'}|)$ belongs to $\cB_1$. This yields the real analyticity of $w$. \cqfd

\subsection{Twisted regular differential calculus.}
\label{ss:twisted-regu}

Next, we let act the twist $U$ by conjugation on the regular, differential part $\uP\in{\rm Diff}_2(\uOmega)$ of Hamiltonians $\uH$ in the class ${\rm Ham}_2(\uOmega)$ (cf. Definitions~\ref{def:ope-diff} and~\ref{def:hamiltonians}). It turns out that this action preserves this class ${\rm Diff}_2(\uOmega)$. 

\begin{proposition}\label{p:twisted-op}
Let $\rx^0\in\cU\cap\uOmega\subset\R^{3m}$. Take $r_0>0$ and $\delta _0\in ]0; r_0/2]$ such that \eqref{eq:écartement} and \eqref{eq:dzf-1} hold true, and such that $\Omega (\delta_0):=B(\rx_1^0; \delta _0[\times\cdots\times B(\rx_m^0; \delta _0[\subset\uOmega$. We recall that $\Omega (\delta_0)\subset\cU$. Let us define $U^\ast : \Omega (\delta_0)\ni\rx\donne U_\rx^{-1}\in\cB_0$. \\
Let $\uP\in{\rm Diff}_2(\uOmega)$ {\rm (}see Definition~\ref{def:ope-diff}{\rm )}.Then the operator composition $P=U\uP U^\ast$ defines a differential operator on $\Omega (\delta_0)\times\R^{3p}$ that belongs to ${\rm Diff}_2(\Omega (\delta_0))$. 
\end{proposition}
\begin{remark}\label{r:ellipticity}
In {\rm\cite{hu,j}}, it was proven that the twisted Laplacian $U\Delta U^\ast$ is elliptic, when $m=1$. Proposition~\ref{p:twisted-op} shows that the conjugation by the twist of an elliptic operator is also elliptic. This property is actually used in {\rm\cite{kmsw,ms}} in a semiclassical framework. 
\end{remark}
\begin{remark}\label{r:oif}
In the proof of Proposition~\ref{p:twisted-op} below, we essentially follow 
{\rm\cite{hu,j}}. Observing that $U$ is a Fourier integral operator {\rm (}cf. {\rm\cite{gs,h4}}{\rm )}, one can use composition rules for such operators to derive the relationship between the principal symbols of $\uP$ and $P$, given in \eqref{eq:symb-princ-twist} below. 
\end{remark}

{\bf Proof of Proposition~\ref{p:twisted-op}:} 
Take $\uP\in{\rm Diff}_2(\uOmega)$, that is of the form \eqref{eq:op-diff2}. 
We first observe that, for $(\alpha ; \, \beta)\in\N^{3m}\times\N^{3p}$ with $|\alpha |+|\beta |\leq 2$, $Uc_{\alpha\beta}U^\ast$ is, as differential operator, the multiplication operator by the function 
\[\Omega (\delta_0)\times\R^{3p}\ni (\rx ; \ry)\ \donne\ c_{\alpha\beta}\bigl(\rx ; F(\rx ; \ry)\bigr)\comma\]
which belongs to $\cC_b^\infty(\Omega (\delta_0)\times\R^{3p}; \C)$.
Using \eqref{eq:twisted-2-gradiant-x} and \eqref{eq:twisted-2-gradiant-y} 
componentwise and taking advantage of \eqref{eq:prop-J-1} and \eqref{eq:prop-J-2} for $G=F^{\langle -1\rangle}$, we see that $P=U\uP U^\ast$ has the form \eqref{eq:op-diff2} on $\Omega (\delta_0)\times\R^{3p}$ with coefficients in $\cC_b^\infty(\Omega (\delta_0)\times\R^{3p}; \C)$. \\
We are left with the proof of the ellipticity of $P$. To this end, we compute its principal symbol $\sigma _P$. Let us denote by ``$\cdot$'' the scalar product in $\R^d$, for any $d\in\N^\ast$. 
It is well-known that the total symbol of $P$ is the map $S_P: \Omega (\delta_0)\times\R^{3p}\times\R^{3m}\times\R^{3p}\dans\C$ defined by 
\begin{equation}\label{eq:symbol}
 S_P(\rx ; \ry; \xi; \eta)\ =\ e^{-i\rx\cdot\xi-i\ry\cdot\eta}P\bigl(e^{i\rx\cdot\xi+i\ry\cdot\eta}\bigr)\period
\end{equation}
Using this, we obtain in Appendix~\ref{app:computations} that, for $(\rx ; \ry; \xi; \eta)\in\Omega (\delta_0)\times\R^{3p}\times\R^{3m}\times\R^{3p}$, 
\begin{equation}\label{eq:symb-princ-twist}
 \sigma _P(\rx ; \ry; \xi; \eta)\ =\ \sigma _\uP \bigl(\rx;\, F(\rx ; \ry); \, \xi + J_1\bigl(F^{\langle -1\rangle}\bigr)(\rx ; \ry)\, \eta ;\,  J_3\bigl(F^{\langle -1\rangle}\bigr)(\rx ; \ry)\, \eta \bigr)\period
\end{equation}
Since $\uP$ is elliptic, there exists $\uC >0$ such that \eqref{eq:glob-ellipt} holds true with $P$ replaced by $\uP$ and $C$ replaced by $\uC$. Take $(\xi ; \eta)$ satisfying $|\xi |^2+|\eta|^2\geq 1$. Then 
\begin{equation}\label{eq:mino-symbol}
 \uC^{-1}\sigma _P(\rx ; \ry; \xi; \eta)\ \geq \ \bigl|\xi\, +\, J_1\bigl(F^{\langle -1\rangle}\bigr)(\rx ; \ry)\eta\bigr|^2\, +\, \bigl|J_3\bigl(F^{\langle -1\rangle}\bigr)(\rx ; \ry)\eta\bigr|^2\period
\end{equation}
By \eqref{eq:prop-J-1}, we can find $M>0$ such that, for all $(\rx ; \ry)\in\Omega (\delta_0)\times\R^{3p}$, 
\[\bigl\|J_1\bigl(F^{\langle -1\rangle}\bigr)(\rx ; \ry)\bigr\|_{\cL(\R^{3p}; \R^{3m})}\ \leq \ M\period\]
Now let $S=\sqrt{1+4M^2}$. Consider first the case where $S|\eta|\leq (|\xi|^2+|\eta |^2)^{1/2}$. We have 
$2M|\eta|\leq\ |\xi|$. Thus 
\[\bigl|\xi\, +\, J_1\bigl(F^{\langle -1\rangle}\bigr)(\rx ; \ry)\eta\bigr|^2\ \geq \ \frac{|\xi|^2}{4}\]
and, using \eqref{eq:inv-J3}, we obtain from \eqref{eq:mino-symbol} the lower bound 
\[\uC^{-1}\sigma _P(\rx ; \ry; \xi; \eta)\geq \min (1/4; C^{-2})(|\xi|^2+|\eta |^2)\period\]
If, now, $S|\eta|\geq (|\xi|^2+|\eta |^2)^{1/2}$, it follows from \eqref{eq:inv-J3} and \eqref{eq:mino-symbol} that 
\[\uC^{-1}\sigma _P(\rx ; \ry; \xi; \eta)\geq C^{-2}S^{-2}(|\xi|^2+|\eta |^2)\period\]
This yields \eqref{eq:glob-ellipt}. Thus $P\in{\rm Diff}_2(\Omega (\delta_0))$. \cqfd

\subsection{Elliptic regularity for twisted Hamiltonians.}
\label{ss:twist-elliptic-regu}

Now, we conjugate by the twist $U$ any Hamiltonians $\uH$ of the class ${\rm Ham}_2(\uOmega)$ (cf. Definition~\ref{def:hamiltonians}) and show that we can apply Theorem~\ref{th:ellipti-regu} to the twisted operator. 
\begin{theorem}\label{th:regu-twisted-op}
Let $\uOmega$ be an open subset of $\R^{3m}$ and $\rx^0\in\cU\cap\uOmega$. Then, one can find a neighbourhood $\Omega _0$ of $\rx^0$ and a unitary operators valued map $U : \Omega _0\ni\rx\donne U_\rx\in\cB_0$ such that the following holds true: \\
For any Hamiltonian $\uH\in {\rm Ham}_2(\uOmega)$ on $\uOmega\times\R^{3p}\subset\R^{3(m+p)}$, for any $\Psi\in W^{2,2}(\uOmega\times\R^{3p})$ such that $\uH\Psi=0$, the map $U\Psi : \Omega _0\ni\rx\donne U_\rx\Psi (\rx ; \cdot)\in\cW_2$ is well-defined and real analytic. Furthermore, for any $\alpha\in\N^{3m}$ with $1\leq |\alpha|\leq 2$, the map $U(D_\rx^\alpha\Psi) : \Omega_0\ni\rx\donne U_\rx (D_\rx^\alpha\Psi) (\rx ; \cdot)\in\cW_{2-|\alpha|}$ is well-defined and real analytic. 
\end{theorem}
{\bf Proof:} Given $\rx^0\in\cU\cap\uOmega$, take $\Omega_0=\Omega (\delta _0)$ satisfying the assumptions of Proposition~\ref{p:twisted-op} and consider the twist $U : \Omega _0\to\cB_0$ defined by \eqref{eq:def-U}. Recall that $U^\ast : \Omega_0\ni\rx\donne U_\rx^{-1}\in\cB_0$. Take $\uH=\uP+\uV\in {\rm Ham}_2(\uOmega)$ and $\Psi\in W^{2,2}(\uOmega\times\R^{3p})$ such that $\uH\Psi=0$. 
Applying Proposition~\ref{p:pot-regu-anal} to each term of $\uV$ and Proposition~\ref{p:twisted-op} to $\uP$, we see that $P=U\uP U^\ast$ and $W=U\uV U^\ast$ satisfy the assumptions of Theorem~\ref{th:ellipti-regu} on $\Omega _0$, yielding the real analyticity of the map $U\Psi : \Omega _0\ni\rx\donne U_\rx\Psi (\rx ; \cdot)\in\cW_2$. \\
By \eqref{eq:twisted-2-gradiant-x}, we have, for $\rx\in\Omega _0$, setting $J_1:=J_1(F^{\langle -1\rangle})$ and $J_2:=J_2(F^{\langle -1\rangle})$, 
\begin{align*}
 U_\rx\bigl(D_\rx\Psi\bigr)(\rx ; \cdot)\ =&\ U_\rx D_\rx U_\rx^{-1}\cdot U_\rx\Psi(\rx ; \cdot)\\ 
 =&\ D_\rx\bigl(U_\rx\Psi(\rx ; \cdot)\bigr)\, +\, J_1(\rx ; \cdot)\cdot D_\ry\bigl(U_\rx\Psi(\rx ; \cdot)\bigr)\, +\, J_2(\rx ; \cdot)\bigl(U_\rx\Psi(\rx ; \cdot)\bigr)\period
\end{align*}
By \eqref{eq:prop-J-1}, this implies that, for $\alpha\in\N^{3m}$ with $|\alpha|=1$, $U(D_\rx^\alpha\Psi) : \Omega_0\ni\rx\donne U_\rx (D_\rx^\alpha\Psi) (\rx ; \cdot)\in\cW_1$ is well-defined and real analytic. Letting act 
$U_\rx D_\rx U_\rx^{-1}$ on the above equality and using \eqref{eq:twisted-2-gradiant-x}, \eqref{eq:prop-J-1}, and \eqref{eq:prop-J-2}, we obtain that, 
for $\alpha\in\N^{3m}$ with $|\alpha|=2$, $U(D_\rx^\alpha\Psi) : \Omega_0\ni\rx\donne U_\rx (D_\rx^\alpha\Psi) (\rx ; \cdot)\in\cW_0$ is well-defined and real analytic. \cqfd

\begin{remark}\label{r:kmsw-hu}
We observe that the present framework {\rm (}the one in Theorem~\ref{th:regu-twisted-op}{\rm )} is closer to the one in {\rm\cite{kmsw}} than the one in Hunziker's original paper {\rm\cite{hu}}. This comes from the fact that, here and in {\rm\cite{kmsw}}, the considered Hamiltonian contains derivatives w.r.t. the ``external'' variable $\rx$, whereas it is not the case in {\rm\cite{hu}}. We further point out that a variant of Hunziker's twist is used in a similar way in {\rm\cite{mm}}. 
\end{remark}
%

\section{Proof of the main results.}
\label{s:proof}
\setcounter{equation}{0}

In this section, we provide a proof for our main results, namely Theorem~\ref{th:courant-anal}, 
Theorem~\ref{th:rho-k-anal}, and Theorem~\ref{th:gamma-k-anal}. We first prove the result on the current density $C$. 

{\bf Proof of Theorem~\ref{th:courant-anal}:} 
Let $x^0\in\R^3\setminus\{R_1,\cdots , R_L\}$. Notice that, if we set $m=1$ and $p=N-1$, then $\cU=\R^3\setminus\{R_1,\cdots , R_L\}$, by \eqref{eq:def-U-ronde}. Note that $\uH:=H-E\in {\rm Ham}(\R^3)$ and $\uH\psi=0$. 
Applying Theorem~\ref{th:regu-twisted-op} to $\uOmega=\R^{3m}=\R^3$ and $\rx^0=x^0\in\cU\cap\uOmega$, there exists some neighbourhood $\Omega_0$ of $x^0$ and a map $U : \Omega _0\ni x\donne U_x\in\cB_0$, the values of which are unitary operators on $\cW_0$, such that, for all $\alpha\in\N^{3}$ with $|\alpha|\leq 1$, the map $U(D_x^\alpha \psi ): \Omega_0\ni x\donne U_x(D_x^\alpha\psi) (x ; \cdot)\in\cW_{2-|\alpha|}$ is well-defined and real analytic. \\
Let $\alpha\in\N^{3}$ with $|\alpha|=1$. For $x\in\Omega _0$, we can write, since $U_x$ is unitary,  
\[\bigl\langle \bigl(D_x^\alpha\psi\bigr) (x; \cdot )\, ,\, \psi (x; \cdot)\bigr\rangle\ =\ \bigl\langle U_x\bigl(D_x^\alpha\psi\bigr) (x; \cdot )\, ,\, U_x\psi (x; \cdot)\bigr\rangle\period\]
As the scalar product of real analytic, $\cW_0$-valued maps, the map 
$\Omega _0\ni x\donne \langle (D_x^\alpha\psi), \psi\rangle$ is also real analytic. This gives the desired result by \eqref{eq:courant-densité}. \cqfd

Next we come to the proof of Theorem~\ref{th:rho-k-anal} on the densities $\rho _k$. 
We also comment on the limitation on the analyticity domain in this result (cf. Remark~\ref{r:limitation-1}). 

{\bf Proof of Theorem~\ref{th:rho-k-anal}:} Take an integer $k$ such that $0<k<N$. Let $\ux ^0=(x_1^0; \cdots ; x_k^0)\in\cU^{(1)}_k$. Setting $m=k$ and $p=N-k$, we observe that $\cU^{(1)}_k$ coincide with the set $\cU$, defined in \eqref{eq:def-U-ronde}.
Since $-\Delta =-\Delta _\ux - \Delta _y\in {\rm Diff}_2(\R^{3m})$ and $V-E$ is a potential of the type $\uV$ in \eqref{eq:def-V} on $\R^{3m}\times\R^{3p}$, by Remark~\ref{r:coulomb-twisté}, we may apply Theorem~\ref{th:regu-twisted-op} to $\uH=H-E$, $\uOmega =\R^{3m}$, $\Psi =\psi$, and $\rx^0=\ux ^0\in\cU\cap\uOmega$. Thus there exist some neighbourhood $\Omega _0$ of $\ux ^0$ and a map $U : \Omega _0\to\cB_0$, with values in the set of unitary operators on $\cW_0$, such that 
the map $U\psi : \Omega_0\ni\ux\donne U_\ux\psi (\ux ;\cdot)\in\cW_2$ is well-defined and real analytic. For $\ux\in\Omega _0$, the density $\rho _k(\ux)$ is, by definition (see \eqref{eq:rho-densité}), the squared $\cW_0$-norm of $\psi (\ux; \cdot)$, which is also the one of 
$U_\ux\psi (\ux; \cdot)$, since $U_\ux$ is unitary on $\cW_0$. Thus, $\rho _k$ is real analytic on $\Omega _0$. \cqfd
\begin{remark}\label{r:basique-1}
 In the $k=1$ case, the above proof reproduces the one of {\rm \cite{j}} on the density $\rho=\rho _1$. Thanks to Hunziker's twist, the proof in the general case is essentially identical.
\end{remark}
\begin{remark}\label{r:limitation-1}
Theorem~\ref{th:rho-k-anal} ensures the real analyticity of $\rho _k$ on $\cU^{(1)}_k$. We actually can give two obstacles to the extension of this property to a larger set by the method used in Section~\ref{s:twisted-diff}.\\ 
If we try to apply this method near a point $\ux ^0\in\cC_k\cup\cR _k$, we have to work near the singularity of a term of the form $|x_j-R_\ell|^{-1}$ or $|x_j-x_{j'}|^{-1}$ with $1\leq j\neq j'\leq k$ and such a singular term would be unaffected by any twist. Thus, we cannot apply Theorem~\ref{th:ellipti-regu} to the twisted equation to get the result of Theorem~\ref{th:regu-twisted-op}.\\
The second obstacle takes place in the construction of the twist near such a point $\ux ^0\in\cC_k\cup\cR _k$. To explain this, let us use the notation of Subsection~\ref{ss:twist} and consider a point $\rx^0$ such that, either $\rx_j^0=R_\ell$ or $\rx_j^0=\rx_{j'}^0$, for 
$j\neq j'$. Assume that we have a function $f$ defined on a neighbourhood of such $\rx ^0$ times $\R^3$ such that, for $\rx$ close to $\rx^0$, $f(\rx ; \cdot)$ is a diffeomorphism on $\R^3$ that satisfies \eqref{eq:transformation}. If $\rx _j^0=\rx _{j'}^0$, then, for $\rx$ close to $\rx ^0$ with different $\rx _j$ and $\rx _{j'}$, we would have $\rx _j=f(\rx; \rx_j^0)=f(\rx; \rx_{j'}^0)=\rx _{j'}$, which is a contradiction. If $\rx _j^0=R_\ell$, then, for $\rx$ close to $\rx ^0$ with $\rx _j\neq R_\ell$, we would have $\rx _j=f(\rx; \rx_j^0)=f(\rx; R_\ell)=R_\ell$, which is again a contradiction.
\end{remark}

Finally, we prove Theorem~\ref{th:gamma-k-anal} concerning the density matrices $\gamma_k$. 
Again we provide a comment on the limitation on the analyticity domain in this result (cf. Remark~\ref{r:limitation-2}).

{\bf Proof of Theorem~\ref{th:gamma-k-anal}:} Let $k$ be an integer between $0$ and $N$. Let $\ux ^0=(x_1^0; \cdots ; x_k^0)$ and $\ux ^{'0}=(x_1^{'0}; \cdots ; x_k^{'0})$ two points in $(\R^3)^k$ such that $(\ux ^0; \ux ^{'0})\in\cU_k^{(2)}$. Setting $m=2k$ and $p=N-k$, we note that $\cU=\cU_k^{(2)}$. 
Denoting by $y$ the variable in $\R^{3(N-k)}=\R^{3p}$, we consider the Hamiltonians $H_1$ and $H_2$ on $(\R^{3k})^2\times\R^{3(N-k)}=\R^{3m}\times\R^{3p}$ defined by
\[H_1\ =\ -\Delta _{\ux}\, -\, \Delta _{\ux '}\, -\, \Delta _{y}\, +\, V(\ux ; y)\hspace{.4cm}\mbox{and}\hspace{0.4cm}H_2\ =\ -\Delta _{\ux}\, -\, \Delta _{\ux '}\, -\, \Delta _{y}\, +\, V(\ux '; y)\]
and the almost everywhere defined functions $\Psi _1 : \R^{3m}\times\R^{3p}\ni (\ux ; \ux '; y)\donne\psi (\ux ; y)$ and $\Psi _2 : \R^{3m}\times\R^{3p}\ni (\ux ; \ux '; y)\donne\psi (\ux '; y)$. 
We note that $H_1\in{\rm Ham}_2(\R^{3m})$ and $H_2\in{\rm Ham}_2(\R^{3m})$. 
From $H\psi =E\psi$, we derive that $H_1\Psi_1=E\Psi_1$ and $H_2\Psi_2=E\Psi_2$, in the distributional sense. Furthermore, for any bounded subset $B$ of $\R^{3m}$, we observe that 
\[\Psi _1\in W^{2,2}\bigl(B\times\R^{3p}_y\bigr)\hspace{.4cm}\mbox{and}\hspace{.4cm}\Psi _2\in W^{2,2}\bigl(B\times\R^{3p}_y\bigr)\period\]
Let $B_0$ be a bounded neighbourhood of $(\ux ^0; \ux ^{'0})$. \\
Now we apply Theorem~\ref{th:regu-twisted-op} to $\uOmega =B_0$ and $\rx^0=(\ux^0 ; \ux^{'0})\in\cU\cap\uOmega$. Therefore there exist a bounded neighbourhood $\Omega _0$ of $(\ux^0 ; \ux^{'0})$ and a map $U : \Omega _0\to\cB_0$, with values in the set of unitary operators on $\cW_0$, such that the maps $U\Psi _1: \Omega_0\ni (\ux ; \ux ')\donne U_{\ux ; \ux '}\Psi _1(\ux ; \ux ';\cdot)\in\cW_2$ and $U\Psi _2: \Omega_0\ni (\ux ; \ux ')\donne U_{\ux ; \ux '}\Psi _2(\ux ; \ux ';\cdot)\in\cW_2$ are well-defined and real analytic.
Now, by definition of $\gamma _k$ (see \eqref{eq:gamma-densité}), we have, for $(\ux ; \ux ')\in\Omega _0$, 
\[\gamma _k(\ux ; \ux ')\ =\ \bigl\langle \Psi _1(\ux ; \ux '; \cdot )\, , \, \Psi _2(\ux ; \ux '; \cdot )\bigr\rangle\ =\ \bigl\langle U_{\ux ; \ux '}\Psi _1(\ux ; \ux '; \cdot )\, , \, U_{\ux ; \ux '}\Psi _2(\ux ; \ux '; \cdot )\bigr\rangle\comma\]
since $U_{\ux ; \ux '}$ is unitary. As the scalar product of real analytic, $\cW_0$-valued maps on $\Omega _0$, $\gamma_k$ is also real analytic on $\Omega _0$. \cqfd

\begin{remark}\label{r:basique-2}
We note that our proofs of Theorems~\ref{th:courant-anal}, \ref{th:rho-k-anal}, and~\ref{th:gamma-k-anal}, have a common structure. They all use an appropriate twist and Theorem~\ref{th:regu-twisted-op}. Differences between them occur in the used set of variables. 
\end{remark}
\begin{remark}\label{r:limitation-2}
According to the first obstacle mentioned in Remark~\ref{r:limitation-1}, the method of Section~\ref{s:twisted-diff} requires to work on $\cU_k^{(1)}\times\cU_k^{(1)}$. Now, because of the second obstacle mentioned there, we also need to exclude the set $\cC _k^{(2)}$ and, therefore, work on $\cU _k^{(2)}=(\cU_k^{(1)}\times\cU_k^{(1)})\setminus\cC _k^{(2)}$. 
\end{remark}
%

\section{Extensions.}
\label{s:extensions}
\setcounter{equation}{0}

In Section~\ref{s:twisted-diff}, we worked in a larger framework than the one that would be needed to treat the physical operator $H$ in \eqref{eq:hamiltonien}. It is thus natural to expect that our main results on $H$ extend to a larger class of Hamiltonians. This is the case, as we shall see in this Section. 

We consider again a system of $N$ electrons, interacting to one another and moving under the influence of $L$ fixed nuclei. Let $\cI :=\{1; \cdots ; N\}$ be the set of indices for the electronic variables $x\in\R^{3N}$ and let $\cI _n:=\{1; \cdots ; L\}$ be the one of the nuclear variables $R_\ell\in\R^3$. It is convenient to introduce the disjoint union of those sets, namely 
\[\cI\, \sqcup\, \cI _n\ :=\ \bigl\{(a; j)\, ;\ a\in\{e; n\}\comma\ j\in\cI_a\bigr\}\period\]
For $c=(a; j)\in \cI _e\sqcup\cI _n$, let $\pi_1(c)=a$ and $\pi_2(c)=j$. We define $X_c=x_j$ if $\pi_1(c)=e$ and $X_c=R_j$ if $\pi_1(c)=n$. 
Let 
\[\widetilde\cP \ :=\ \bigl\{\{c; c'\}\, ;\ (c; c')\in \bigl(\cI _e\sqcup\, \cI _n\bigr)^2\ \mbox{and}\ c\neq c'\bigr\}\]
be the set of all possible particle pairings. We consider a potential $\tV$ where 
\begin{equation}\label{eq:def-tilde-V}
\forall x\in\R^{3N}\comma\hspace{.4cm}\tV (x)\ =\ \sum _{\{c;\, c'\}\in\widetilde\cP}\, V_{\{c;\, c'\}}(X_c\, -\, X_{c'})\comma
\end{equation}
with $V_{\{c;\, c'\}}\in\cV$ (cf. Definition~\ref{def:classe-V}), for all $\{c; c'\}\in\widetilde\cP$. \\
Now, let $\tP=\tP (x; D_x)$ a differential operator on $\R^{3N}$ of the form 
\begin{equation}\label{eq:op-tilde-diff2}
 \tP\ =\ \sum_{\stackrel{\alpha \in\N^{3N}}{\stackrel{}{|\alpha |\leq 2}}}\, \tc_{\alpha}(x)\, D_x^\alpha \comma
\end{equation}
where the coefficients $\tc_{\alpha}$ belong to $\cC_b^\infty(\R^{3N}; \C)$. Let $\tH=\tP+\tV$ and assume that there exists $\tpsi\in W^{2,2}(\R^{3N})$ such that $\tH\tpsi=0$. For $0<k<N$, we still denote by $\rho _k$, $\gamma _k$, and $C$ the objects that are respectively defined by \eqref{eq:rho-densité}, \eqref{eq:gamma-densité}, and \eqref{eq:courant-densité}, 
with $\psi$ replaced by $\tpsi$. 
\begin{theorem}\label{th:géné-anal}
Under the above assumptions on $\tH$ and $\tpsi$, the results of Theorems~\ref{th:courant-anal}, \ref{th:rho-k-anal}, and~\ref{th:gamma-k-anal}, are valid. 
\end{theorem}
{\bf Proof:} We can follow the above proofs of Theorems~\ref{th:courant-anal}, \ref{th:rho-k-anal}, and~\ref{th:gamma-k-anal}. In the later, we use the elliptic operators 
\[H_1\ =\ -\Delta _{\ux '}\, -\, \tP (\ux; y; D_\ux; D_y)\, +\, \tV(\ux ; y)\hspace{.4cm}\mbox{and}\hspace{0.4cm}H_2\ =\ -\Delta _\ux\, -\, \tP (\ux '; y; D_{\ux '}; D_y)\, +\, \tV(\ux '; y)\]
and the almost everywhere defined functions $\Psi _1 : \R^{3m}\times\R^{3p}\ni (\ux ; \ux '; y)\donne\tpsi (\ux ; y)$ and $\Psi _2 : \R^{3m}\times\R^{3p}\ni (\ux ; \ux '; y)\donne\tpsi (\ux '; y)$. \cqfd

\begin{remark}\label{r:extension}
We point out that the present class of operators includes Laplace-Beltrami operators associated to appropriate metrics. One can also replace $D_x$ by $D_x-A(x)$, where $A$ is a real analytic vector potential in $\cC_b^\infty(\R^{3N}; \R^3)$, and add an external electric potential $V_e$ in $\cC_b^\infty(\R^{3N}; \C)$, that is also real analytic. \\
Notice that the class $\cV$ in Definition~\ref{def:classe-V} contains non-radial pair potentials. Among the admissible functions $\eta$ used to define $\cV$, we have $0<t\donne t^{-1}$ but we can also choose $0<t\donne t^{-1}(\ln (t))^{-\epsilon}$, with $\epsilon>0$. We observe that this class $\cV$ is contained in the class of pair potentials used in {\rm\cite{fs}}. 
\end{remark}
%


\newpage
\begin{appendices}{\bf \Large Appendix.}

\renewcommand{\theequation}{\thesection .\arabic{equation}}

\appendixtitleon
\appendixtitletocon

\section{Basic computations.}
\label{app:computations}
\setcounter{equation}{0}

In this section, we provide explicit formulae for several objects related to the twist $U$ in Subsection~\ref{ss:twist} and the computation of the principal symbol of some differential operator, that is needed in Subsection~\ref{ss:twisted-regu}. 

First of all, let us notice that, from the expressions of $d_\rx f$ and $d_z f$, there exists 
some $\tau$, $r_0$, and $\rx^0$ dependent constant $M$, such that, denoting by $I_3$ the identity on $\R^3$, 
\begin{equation}\label{eq:support}
\forall\rx\in\R^{3m}\comma\, \forall z\in\R^3\comma\hspace{.4cm}|z|\geq M\ \impl\ \bigl((d_\rx f)(\rx; z)\ =\ 0\hspace{.1cm}\mbox{and}\hspace{.1cm}(d_z f)(\rx; z)\, -\, I_3\ =\ 0\bigr)\period
\end{equation}
Using \eqref{eq:dzf-1}, we can find some $c>0$ such that
\begin{equation}\label{eq:bornes-dzf}
\forall (\rx ; z)\in\Omega (\delta _0)\times\R^{3}\comma\ c^{-1}\, \leq\ \bigl\|\bigl(\bigl(d_zf\bigr)\bigl(\rx ; z\bigr)\bigr)^{-1}\bigr\|_{\cL (\R^3)}\, +\, \bigl\|\bigl(d_zf\bigr)\bigl(\rx ; z\bigr)\bigr\|_{\cL (\R^3)}\ \leq\ c\period
\end{equation}
For $\rx\in\Omega (\delta _0)$, $j\in\ncg 1; m\ncd$, and $z\in\R^3$, we have 
\[f^{\langle -1\rangle}\bigl(\rx ; f(\rx ; z)\bigr)\ =\ z\period\]
By differentiation of this equality and by using \eqref{eq:dxf}, we obtain, with equalities in $\cL(\R^3)$, in $\cL(\R^3)$, and in $\cL (\R^{3m}; \R^3)$, respectively, 
\begin{align}
\bigl(d_zf^{\langle -1\rangle}\bigr)\bigl(\rx ; f(\rx ; z)\bigr)\ &=\ \bigl((d_zf)(\rx ; z)\bigr)^{-1}\comma\label{eq:diff-z-inverse-f}\\
\bigl(d_{\rx_j} f^{\langle -1\rangle}\bigr)\bigl(\rx ; f(\rx ; z)\bigr)\ &=\ -\, \tau\bigl(r_0^{-1}(z\, -\rx_j^0)\bigr)\bigl(d_zf^{\langle -1\rangle}\bigr)\bigl(\rx ; f(\rx ; z)\bigr)\comma\label{eq:diff-xj-inverse-f}\\
\bigl(d_{\rx} f^{\langle -1\rangle}\bigr)\bigl(\rx ; f(\rx ; z)\bigr)\ &=\ -\, \bigl(d_zf^{\langle -1\rangle}\bigr)\bigl(\rx ; f(\rx ; z)\bigr)\circ\sum_{j'=1}^m\tau\bigl(r_0^{-1}(z\, -\rx_{j'}^0)\bigr)\, d\rx _{j'}\period\label{eq:diff-x-inverse-f}
\end{align}
Now, using \eqref{eq:support}, \eqref{eq:bornes-dzf}, \eqref{eq:diff-z-inverse-f}, \eqref{eq:diff-xj-inverse-f}, and \eqref{eq:diff-x-inverse-f}, it is straightforward to check that 
\begin{equation}\label{eq:c-infini-b-f}
 f\in\cC_b^\infty \bigl(\Omega (\delta _0)\times\R^3; \R^3\bigr)\hspace{.4cm}\mbox{and}\hspace{.4cm}f^{\langle -1\rangle}\in\cC_b^\infty \bigl(\Omega (\delta _0)\times\R^3; \R^3\bigr)\period
\end{equation}
From the very definition of $F$, we note that, for $\rx\in\Omega (\delta _0)$, $\rx^\epsilon\in\R^{3m}$, and $\ry\in\R^{3p}$, 
\[(d_\rx F)(\rx ; \ry)\cdot\rx ^\epsilon\ =\ \bigl((d_\rx f)(\rx ; \ry _1)\cdot\rx ^\epsilon \, ;\ \cdots \ ;\,  (d_\rx f)(\rx ; \ry_p)\cdot\rx ^\epsilon\bigr)\, \in\, \R^{3p}\]
and the $k$-th component (with $k\in\ncg 1; p\ncd$) in the latter expression is given by 
\[(d_\rx f)(\rx ; \ry_k)\cdot\rx ^\epsilon\ =\ \sum _{j=1}^m\tau\bigl(r_0^{-1}(\ry _k\, -\rx_j^0)\bigr)\, \rx _j^\epsilon\, \in\, \R^3\comma\]
by \eqref{eq:dxf}. For $\rx\in\Omega (\delta _0)$, $\ry\in\R^{3p}$, $z\in\R^3$, and $k\in\ncg 1; p\ncd$, we have 
\[(d_{\ry_k} F)(\rx ; \ry)\cdot z\ =\ \bigl(0\, ;\, \cdots \, ;\, 0\, ;\, (d_z f)(\rx ; \ry _k)\cdot z\, ;\, 0 \, ;\ \cdots \ ;\,  0\bigr)\, \in\, \R^{3p}\]
(the term containing $z$ being in the $k$th position) and, for $\ry^\epsilon\in\R^{3p}$, 
\[ (d_{\ry} F)(\rx ; \ry)\cdot\ry^\epsilon\ =\ \sum _{k=1}^p\, \bigl(0\, ;\, \cdots \, ;\, 0\, ;\, (d_z f)(\rx ; \ry _k)\, \cdot \, \ry^\epsilon_k\, ;\, 0 \, ;\ \cdots \ ;\,  0\bigr)\, \in\, \R^{3p}\]
\[=\ \ry^\epsilon \, +\, r_0^{-1}\sum _{k=1}^p\, \sum _{j=1}^m\Bigl(0\, ;\, \cdots \, ;\, 0\, ;\, \bigl((d\tau) \bigl(r_0^{-1}(\ry_k\, -\rx_j^0)\bigr)\cdot \ry^\epsilon _k\bigr)(\rx _j\, -\, \rx _j^0)\, ;\, 0 \, ;\ \cdots \ ;\,  0\Bigr)\comma\]
by \eqref{eq:dzf}, each term containing $d\tau$ being in the $k$th position. 
Since we can write, for $\rx\in\Omega (\delta _0)$ and $\ry\in\R^{3p}$, 
\[F^{\langle -1\rangle}\bigl(\rx ; \ry\bigr)\ =\ \bigl(f^{\langle -1\rangle}(\rx ; \ry _1)\, ;\, \cdots\, ;\, f^{\langle -1\rangle}(\rx ; \ry _p)\bigr)\, \in \, \R^{3p}\comma\]
we have, for $\rx^\epsilon\in\R^{3m}$, 
\begin{align*}
 &\bigl(d_\rx F^{\langle -1\rangle}\bigr)\bigl(\rx ; F(\rx ; \ry)\bigr)\, \cdot\, \rx^\epsilon\\
 =&\ \bigl(\bigl(d_\rx f^{\langle -1\rangle}\bigr)\bigl(\rx ; f(\rx ; \ry _1)\bigr)\cdot\rx ^\epsilon \, ;\ \cdots \ ;\,  \bigl(d_\rx f^{\langle -1\rangle}\bigr)\bigl(\rx ; f(\rx ; \ry _p)\bigr)\cdot\rx ^\epsilon\bigr)\, \in\, \R^{3p}\comma
\end{align*}
where $d_\rx f^{\langle -1\rangle}$ is given by \eqref{eq:diff-x-inverse-f}.
Finally, for $\ry^\epsilon\in\R^{3p}$, we have 
\begin{align*}
 &\bigl(d_\ry F^{\langle -1\rangle}\bigr)\bigl(\rx ; F(\rx ; \ry)\bigr)\, \cdot\, \ry^\epsilon\\
 =&\ \sum _{k=1}^p\, \bigl(0\, ;\, \cdots \, ;\, 0\, ;\, \bigl(d_z f^{\langle -1\rangle}\bigr)\bigl(\rx ; f(\rx ; \ry _k)\bigr)\, \cdot \, \ry^\epsilon_k\, ;\, 0 \, ;\ \cdots \ ;\,  0\bigr)\, \in\, \R^{3p}\comma
\end{align*}
where $d_z f^{\langle -1\rangle}$ is explicited in \eqref{eq:diff-z-inverse-f}.  \\
We derive from \eqref{eq:c-infini-b-f} that 
\begin{equation}\label{eq:c-infini-b-F}
 F\in\cC_b^\infty \bigl(\Omega (\delta _0)\times\R^{3p}; \R^{3p}\bigr)\hspace{.4cm}\mbox{and}\hspace{.4cm}F^{\langle -1\rangle}\in\cC_b^\infty \bigl(\Omega (\delta _0)\times\R^{3p}; \R^{3p}\bigr)\period
\end{equation}

Next, we compute the terms $J_k(G)$, for $1\leq k\leq 4$ and $G\in\{F; F^{\langle -1\rangle}\}$, that appear in \eqref{eq:twisted-1-gradiant-x}, \eqref{eq:twisted-2-gradiant-x}, \eqref{eq:twisted-1-gradiant-y}, and \eqref{eq:twisted-2-gradiant-y}. \\
Let $G\in\{F; F^{\langle -1\rangle}\}$. Let us consider the smooth functions $\rho_\pm\in\cC^\infty (\Omega (\delta _0)\times\R^{3p}; \R)$ defined by 
\[\rho _+(\rx ; \ry)\ =\ \bigl|{\rm Det}\, (d_yG)(\rx ; \ry)\bigr|^{1/2}\hspace{.4cm}\mbox{and}\hspace{.4cm}\rho _-(\rx ; \ry)\ =\ \bigl|{\rm Det}\, \bigl(d_yG^{\langle -1\rangle}\bigr)(\rx ; \ry)\bigr|^{1/2}\period\]
From the identity $G^{\langle -1\rangle}(\rx ; G(\rx ; \ry))=(\rx ; \ry)$, we derive that 
$\rho _+(\rx ; \ry)\rho _-(\rx ; G(\rx ; \ry))=1$. Let $\varphi\in\cC^\infty (\Omega (\delta _0)\times\R^{3p}; \R)$. We set 
\[(U_G\varphi )(\rx ; \ry)\ =\ \rho _+(\rx ; \ry)\, \varphi\bigl(\rx ; G(\rx ; \ry)\bigr)\hspace{.4cm}\mbox{and}\hspace{.4cm}\bigl(U_G^{\langle -1\rangle}\varphi \bigr)(\rx ; \ry)\ =\ \rho _-(\rx ; \ry)\, \varphi\bigl(\rx ; G^{\langle -1\rangle}(\rx ; \ry)\bigr)\period\]
Since 
\begin{align*}
 d_\rx\bigl(U_G^{\langle -1\rangle}\varphi \bigr)(\rx ; \ry)\ =&\ \varphi\bigl(\rx ; G^{\langle -1\rangle}(\rx ; \ry)\bigr)\, d_\rx\rho_-(\rx ; \ry)\\
 &\ +\ \rho_-(\rx ; \ry)\, (d_\rx\varphi )\bigl(\rx ; G^{\langle -1\rangle}(\rx ; \ry)\bigr)\\
 &\ +\ \rho_-(\rx ; \ry)\, (d_\ry\varphi )\bigl(\rx ; G^{\langle -1\rangle}(\rx ; \ry)\bigr)\bigl(d_\rx G^{\langle -1\rangle}\bigr)(\rx ; \ry)\comma
\end{align*}
we have, for $\rx ^\epsilon\in\R^{3m}$, 
\begin{align*}
 d_\rx\bigl(U_G^{\langle -1\rangle}\varphi \bigr)(\rx ; \ry)\rx^\epsilon\ =&\ \varphi\bigl(\rx ; G^{\langle -1\rangle}(\rx ; \ry)\bigr)\, \nabla_\rx\rho_-(\rx ; \ry)\cdot\rx ^\epsilon\\
 &\ +\ \rho_-(\rx ; \ry)\, (\nabla_\rx\varphi )\bigl(\rx ; G^{\langle -1\rangle}(\rx ; \ry)\bigr)\cdot\rx ^\epsilon\\
 &\ +\ \rho_-(\rx ; \ry)\, (\nabla_\ry\varphi )\bigl(\rx ; G^{\langle -1\rangle}(\rx ; \ry)\bigr)\cdot\bigl(d_\rx G^{\langle -1\rangle}\bigr)(\rx ; \ry)\rx ^\epsilon\comma
\end{align*}
with scalar products in $\R^{3m}$, $\R^{3m}$, and $\R^{3p}$, respectively. This yields 
\begin{align*}
 \nabla_\rx\bigl(U_G^{\langle -1\rangle}\varphi \bigr)(\rx ; \ry)\ =&\ \varphi\bigl(\rx ; G^{\langle -1\rangle}(\rx ; \ry)\bigr)\, \nabla_\rx\rho_-(\rx ; \ry)\\
 &\ +\ \rho_-(\rx ; \ry)\, (\nabla_\rx\varphi )\bigl(\rx ; G^{\langle -1\rangle}(\rx ; \ry)\bigr)\\
 &\ +\ \rho_-(\rx ; \ry)\, \bigl(d_\rx G^{\langle -1\rangle}\bigr)^T(\rx ; \ry)(\nabla_\ry\varphi )\bigl(\rx ; G^{\langle -1\rangle}(\rx ; \ry)\bigr)\period
\end{align*}
Thus, using $\rho _+(\rx ; \ry)\rho _-(\rx ; G(\rx ; \ry))=1$, 
\begin{align*}
 \bigl(U_G\nabla_\rx U_G^{\langle -1\rangle}\varphi \bigr)(\rx ; \ry)\ =&\ \rho _+(\rx ; \ry)\varphi (\rx ; \ry)\, (\nabla_\rx\rho_-)\bigl(\rx ; G(\rx ; \ry)\bigr)\ +\ (\nabla_\rx\varphi )(\rx ; \ry)\\
 &\ +\ \bigl(d_\rx G^{\langle -1\rangle}\bigr)^T\bigl(\rx ; G(\rx ; \ry)\bigr)(\nabla_\ry\varphi )(\rx ; \ry)\period
\end{align*}
Since 
\begin{align*}
 d_\ry\bigl(U_G^{\langle -1\rangle}\varphi \bigr)(\rx ; \ry)\ =&\ \varphi\bigl(\rx ; G^{\langle -1\rangle}(\rx ; \ry)\bigr)\, d_\ry\rho_-(\rx ; \ry)\\
 &\ +\ \rho_-(\rx ; \ry)\, (d_\ry\varphi )\bigl(\rx ; G^{\langle -1\rangle}(\rx ; \ry)\bigr)\bigl(d_\ry G^{\langle -1\rangle}\bigr)(\rx ; \ry)\comma
\end{align*}
we have, for $\ry ^\epsilon\in\R^{3p}$, 
\begin{align*}
 d_\ry\bigl(U_G^{\langle -1\rangle}\varphi \bigr)(\rx ; \ry)\ry^\epsilon\ =&\ \varphi\bigl(\rx ; G^{\langle -1\rangle}(\rx ; \ry)\bigr)\, \nabla_\ry\rho_-(\rx ; \ry)\cdot\ry ^\epsilon\\
 &\ +\ \rho_-(\rx ; \ry)\, (\nabla_\ry\varphi )\bigl(\rx ; G^{\langle -1\rangle}(\rx ; \ry)\bigr)\cdot\bigl(d_\ry G^{\langle -1\rangle}\bigr)(\rx ; \ry)\ry ^\epsilon\comma
\end{align*}
with scalar products in $\R^{3p}$. This yields 
\begin{align*}
 \nabla_\ry\bigl(U_G^{\langle -1\rangle}\varphi \bigr)(\rx ; \ry)\ =&\ \varphi\bigl(\rx ; G^{\langle -1\rangle}(\rx ; \ry)\bigr)\, \nabla_\ry\rho_-(\rx ; \ry)\\
 &\ +\ \rho_-(\rx ; \ry)\, \bigl(d_\ry G^{\langle -1\rangle}\bigr)^T(\rx ; \ry)(\nabla_\ry\varphi )\bigl(\rx ; G^{\langle -1\rangle}(\rx ; \ry)\bigr)\period
\end{align*}
Thus, using $\rho _+(\rx ; \ry)\rho _-(\rx ; G(\rx ; \ry))=1$, 
\begin{align*}
 \bigl(U_G\nabla_\ry U_G^{\langle -1\rangle}\varphi \bigr)(\rx ; \ry)\ =&\ \rho _+(\rx ; \ry)\varphi (\rx ; \ry)\, (\nabla_\ry\rho_-)\bigl(\rx ; G(\rx ; \ry)\bigr)\\
 &\ +\ \bigl(d_\rx G^{\langle -1\rangle}\bigr)^T\bigl(\rx ; G(\rx ; \ry)\bigr)(\nabla_\ry\varphi )(\rx ; \ry)\period
\end{align*}
This proves \eqref{eq:twisted-1-gradiant-x}, \eqref{eq:twisted-2-gradiant-x}, \eqref{eq:twisted-1-gradiant-y}, and \eqref{eq:twisted-2-gradiant-y}. Now, taking into account \eqref{eq:c-infini-b-F}, we obtain \eqref{eq:prop-J-1} and~\eqref{eq:prop-J-2}. Using \eqref{eq:bornes-dzf} and~\eqref{eq:diff-z-inverse-f}, we see that \eqref{eq:inv-J3} holds true.

Finally, we compute the principal symbol $\sigma _P$ appearing in the proof of Proposition~\ref{p:twisted-op}. 
In the vector space of polynomials in $\R^{3m}\times\R^{3p}$ with $\cC_b^\infty(\Omega (\delta_0)\times\R^{3p}; \C)$-coefficients, we denote by 
$\equiv$ the equality of polynomials modulo polynomials of degree less than $2$. Then, for $(\rx ; \ry; \xi; \eta)\in\Omega (\delta_0)\times\R^{3p}\times\R^{3m}\times\R^{3p}$, we have, by \eqref{eq:symbol}, \eqref{eq:twisted-2-gradiant-x}, and \eqref{eq:twisted-2-gradiant-y}, 
\begin{align*}
&\sigma _P(\rx ; \ry; \xi; \eta)\\ 
\equiv&\ S_P(\rx ; \ry; \xi; \eta)\\
\equiv&\sum_{\stackrel{(\alpha ; \, \beta)\in\N^{3m}\times\N^{3p}}{\stackrel{}{|\alpha |+|\beta |=2}}}\, c_{\alpha\beta}\bigl(\rx ; F(\rx ; \ry)\bigr)\, \Bigl(\xi\, +\, J_1\bigl(F^{\langle -1\rangle}\bigr)(\rx ; \ry)\, \eta\, +\, J_2\bigl(F^{\langle -1\rangle}\bigr)(\rx ; \ry)\Bigr)^\alpha\\
&\hspace{3.4cm}\times\Bigl(J_3\bigl(F^{\langle -1\rangle}\bigr)(\rx ; \ry)\, \eta\, +\, J_4\bigl(F^{\langle -1\rangle}\bigr)(\rx ; \ry)\Bigr)^\beta\\
&\\
\equiv&\sum_{\stackrel{(\alpha ; \, \beta)\in\N^{3m}\times\N^{3p}}{\stackrel{}{|\alpha |+|\beta |=2}}}\, c_{\alpha\beta}\bigl(\rx ; F(\rx ; \ry)\bigr)\, \Bigl(\xi\, +\, J_1\bigl(F^{\langle -1\rangle}\bigr)(\rx ; \ry)\, \eta\Bigr)^\alpha\Bigl(J_3\bigl(F^{\langle -1\rangle}\bigr)(\rx ; \ry)\, \eta\Bigr)^\beta
\end{align*}
and the last term is actually $\sigma _P(\rx ; \ry; \xi; \eta)$, since it is homogeneous of degree $2$. Thus 
\[\sigma _P(\rx ; \ry; \xi; \eta)\ =\ \sigma _\uP \Bigl(\rx;\, F(\rx ; \ry);\, \xi + J_1\bigl(F^{\langle -1\rangle}\bigr)(\rx ; \ry)\,\eta ;\,  J_3\bigl(F^{\langle -1\rangle}\bigr)(\rx ; \ry)\, \eta \Bigr)\comma\]
yielding \eqref{eq:symb-princ-twist}. 

\section{A global pseudo-differential calculus.}
\label{app:pseudo}
\setcounter{equation}{0}

In this section, we describe a global pseudo-differential calculus that is suited for Section~\ref{s:elliptic}. It is the one developed in the beginning of Chapter 18 in \cite{h2}. 

With the notation of Subsection~\ref{ss:elliptic-regu}, we set $N=m+p$. 
For $k\in\Z$, the class $S^k$ in \cite{h2} (p. 65-75) is the set of smooth, complex valued functions $\sigma$ on $\R^{6N}$ such that, for all $(\alpha ; \beta )\in (\N^{3N})^2$, there exists a constant $C_{\alpha ; \beta}>0$ such that, for all $(\rx ; \ry; \xi ; \eta)\in\R^{6N}$, 
\begin{equation*}
 (1+|\xi |^2+|\eta |^2)^{|\beta |/2}\bigl|\partial _{\rx ; \ry}^\alpha\partial 
_{\xi ; \eta}^{\beta}\sigma(\rx ; \ry; \xi ; \eta)\bigr|\ \leq \ C_{\alpha ; \beta}(1+|\xi |^2+|\eta |^2)^{k/2}\comma
\end{equation*}
where $\partial _{\rx ; \ry}^\alpha$ (resp. $\partial _{\xi ; \eta}^{\beta}$) stands for the $\alpha$ (resp. $\beta$) partial derivative w.r.t. $(\rx ; \ry)$ (resp. $(\xi ; \eta)$). Any such symbol $\sigma$ defines an operator $\sigma (\rx ; \ry; D_\rx ; D_\ry)$ from the Schwarz space $\cS$ on $\R^{3N}$ into itself in the following linear way: for any $u\in\cS$, for any $(\rx ; \ry)\in\R^{3N}$, 
\[\bigl(\sigma (\rx ; \ry; D_\rx ; D_\ry)u\bigr)(\rx ; \ry)\ =\ (2\pi)^{-3N}\int_{\R^{3N}}\, e^{i(\rx ; \ry)\cdot(\xi ; \eta)}\, \sigma(\rx ; \ry; \xi ; \eta)\, \hu (\xi ; \eta)\, d\xi\, d\eta\comma\]
where $(\rx ; \ry)\cdot(\xi ; \eta)$ is the usual scalar product in $\R^{3N}$ of the vectors $(\rx ; \ry)$ and $(\xi ; \eta)$ and where $\hu$ denotes the Fourier transform of $u$. It turns out that, for any $s\in\R$, this operator $\sigma (\rx ; \ry; D_\rx ; D_\ry)$ extends to a bounded operator from $W^{s, 2}(\R^{3N})$ to $W^{s-k, 2}(\R^{3N})$. If $\sigma$ is the constant function equal to one on $\R^{3N}$ then $\sigma\in S^0$ and $\sigma (\rx ; \ry; D_\rx ; D_\ry)$ is the identity operator $I$. \\
For $(k_1; k_2)\in\Z^2$, for $\sigma _1\in S^{k_1}$ and $\sigma _2\in S^{k_2}$, the composition $\sigma _1(\rx ; \ry; D_\rx ; D_\ry)\sigma _2(\rx ; \ry; D_\rx ; D_\ry)$ is given by $\sigma (\rx ; \ry; D_\rx ; D_\ry)$ where $\sigma\in S^{k_1+k_2}$ and $\sigma =\sigma _1\sigma _2+r$ with $r\in S^{k_1+k_2-1}$. 

Now, we are in a position to construct the parametrix $Q$ needed in Subsection~\ref{ss:elliptic-regu}. We observe that $\hP=\sigma (\rx ; \ry; D_\rx ; D_\ry)$ where its symbol $\sigma$ belongs to $S^2$. Furthermore, one can write $\sigma =\sigma_P+r$ with $\sigma _P\in S^2$ and $r\in S^1$. Let $\tau _P\in\cC^\infty _c(\R^{3N})$ such that $\tau (\xi ; \eta)=1$ if $|\xi |^2+|\eta |^2\leq 1$. By the ellipticity assumption of $\hP$ (see \eqref{eq:glob-ellipt}), the symbol 
\[q_1\ :=\ \bigl(1\, -\, \tau _P\bigr)\cdot \sigma _P^{-1}\, \in\, S^{-2}\]
and, setting $Q_1=q_1(\rx ; \ry; D_\rx ; D_\ry)$, there exist symbols $r_1, r_2, r_3\in S^{-1}$ such that 
\begin{align*}
 Q_1\hP\ =\ &(\sigma q_1)(\rx ; \ry; D_\rx ; D_\ry)\, +\, r_1(\rx ; \ry; D_\rx ; D_\ry)\\
 \ =\ &(\sigma _Pq_1)(\rx ; \ry; D_\rx ; D_\ry)\, +\, r_2(\rx ; \ry; D_\rx ; D_\ry)\\
 \ =\ &I\, +\, r_3(\rx ; \ry; D_\rx ; D_\ry)\period
\end{align*}
Let $Q=(I-r_3(\rx ; \ry; D_\rx ; D_\ry))Q_1$. Then, $Q=q(\rx ; \ry; D_\rx ; D_\ry)$ with $q\in S^{-2}$ and there exists 
$r_4\in S^{-2}$ such that $Q\hP=I-r_4(\rx ; \ry; D_\rx ; D_\ry)$. Setting $R:=r_4(\rx ; \ry; D_\rx ; D_\ry)$, we have 
$Q\hP=1-R$ and \eqref{eq:regularisation} holds true. 

\end{appendices}


%
%
\end{document}